\documentclass[12pt,twosides]{amsart}
\usepackage{amssymb,amsmath,amsthm, amscd, enumerate, mathrsfs}
\usepackage{graphicx, hhline}
\usepackage[all]{xy}
\usepackage{braket}
\usepackage[usenames]{color}
\usepackage{hyperref}
\usepackage{fancyhdr}
\usepackage{url}
\usepackage[top=30truemm,bottom=30truemm,left=30truemm,right=30truemm]{geometry}
\def\C{\mathbb C}
\def\R{\mathbb R}
\def\Q{\mathbb Q}
\def\Z{\mathbb Z}

\def\O{\mathcal{O}}

\hypersetup{colorlinks=true}

\title[MMP for lc pairs on complex analytic spaces]{On termination of minimal model program for log canonical pairs on complex analytic spaces}
\author{Makoto Enokizono and Kenta Hashizume}
\date{2025/11/26}
\address{Graduate School of Mathematical Sciences, University of Tokyo,
3-8-1 Komaba, Meguro-ku, Tokyo 153-8914, Japan }
\email{enokizono@g.ecc.u-tokyo.ac.jp}
\address{Department of 
Mathematics, Faculty of Science, Niigata University, Niigata 950-2181, Japan}
\address{Institute for Research Administration, Niigata University, Niigata 950-2181, Japan}
\email{hkenta@math.sc.niigata-u.ac.jp}

\newtheorem{thm}{Theorem}[section]

\newtheorem{lem}[thm]{Lemma}
\newtheorem{cor}[thm]{Corollary}

\newtheorem{conj}[thm]{Conjecture}

\theoremstyle{definition}
\newtheorem{defn}[thm]{Definition}
\newtheorem{rem}[thm]{Remark}

\newtheorem*{ack}{Acknowledgments} 
 
\newtheorem*{b-divisor}{b-divisors} 
 
\newtheorem*{g-pair}{Generalized pairs} 
\newtheorem*{adj-g-pair}{Divisorial adjunction for generalized pairs} 
\newtheorem*{mmp-g-pair}{MMP for generalized pairs}

\newtheorem*{claim*}{Claim}
\begin{document}

\begin{abstract}
We study the termination of minimal model programs for log canonical pairs in the complex analytic setting. 
By using the termination, we prove a relation between the minimal model theory for projective log canonical pairs and that for log canonical pairs in the complex analytic setting. 
The minimal model programs for algebraic stacks and analytic stacks are also discussed.  
\end{abstract}

\maketitle

\tableofcontents

\section{Introduction}

In this paper, we study the termination of a minimal model program (MMP, for short) for log canonical (lc, for short) pairs in the complex analytic setting. 
The finiteness of B-pluricanonical representation (\cite{fujino-gongyo}, \cite{haconxu}) and the existence of lc flips (\cite{birkar-flip}, \cite{haconxu-lcc}) are currently known by Fujino \cite{fujino-analytic-lcabundance}. 
Therefore, we may run an MMP for lc pairs with projective morphisms between complex analytic spaces. 
In the algebraic case, the second author and Hu \cite{hashizumehu} proved the termination of some MMP with scaling of ample divisors for lc pairs when the lc pairs have log minimal models or the log canonical divisors of the lc pairs are not pseudo-effective. 

In this paper, we study the termination of MMP with scaling of ample divisors for projective morphisms between complex analytic spaces.  
The goal of this paper is to prove the complex analytic analog of \cite[Theorem 1.7]{hashizumehu} and to apply the main results to MMP for algebraic stacks and analytic stacks. 
The following theorems are the main results of this paper. 
For the definition of property (P) appearing in the statements below,
see Definition \ref{defn--property(P)}. 

\begin{thm}[cf.~{\cite[Theorem 1.7]{hashizumehu}}]\label{thm--strictmmp-intro} Let $\pi \colon X \to Y$ be a projective morphism from a normal analytic variety $X$ to a Stein space $Y$, and let $W \subset Y$ be a compact subset such that $\pi$ and $W$ satisfy (P). 
Let $(X,\Delta)$ be an lc pair. 
Suppose that $K_{X}+\Delta$ is not pseudo-effective over $Y$ or $(X,\Delta)$ has a log minimal model over $Y$ around $W$ after shrinking $Y$ around $W$. 
Let $A$ be a $\pi$-ample $\mathbb{R}$-divisor on $X$ such that $(X,\Delta+A)$ is lc and $K_{X}+\Delta+A$ is nef over $W$. 
Then, there exist a Stein open subset $Y' \subset Y$ containing $W$ and a sequence of steps of a $(K_{X}+\Delta)$-MMP over $Y'$ around $W$ with scaling of $A$
$$(X,\Delta)=:(X_{0},\Delta_{0})\dashrightarrow (X_{1},\Delta_{1})\dashrightarrow \cdots \dashrightarrow (X_{n},\Delta_{n}),$$
which is represented by bimeromorphic contractions over $Y'$, such that $(X_{n},\Delta_{n})$ is a log minimal model or a Mori fiber space of $(X,\Delta)$ over $Y'$ around $W$. 
Furthermore, if $(X,\Delta)$ has a good minimal model over $Y$ around $W$ after shrinking $Y$ around $W$, then the resulting lc pair $(X_{n},\Delta_{n})$ is a good minimal model of $(X,\Delta)$ over $Y'$ around $W$. 
\end{thm}

\begin{thm}[cf.~{\cite[Theorem 1.1]{birkar-flip}, \cite[Theorem 1.6]{haconxu-lcc}}]\label{thm--mmp-birkar-intro} Let $\pi \colon X \to Y$ be a projective morphism from a normal analytic variety $X$ to a Stein space $Y$, and let $W \subset Y$ be a compact subset such that $\pi$ and $W$ satisfy (P). 
Let $(X,\Delta)$ be an lc pair, and let $B$ be an effective $\mathbb{R}$-Cartier divisor on $X$ such that $K_{X}+\Delta+B\sim_{\mathbb{R},\,Y}0$ and $(X,\Delta+tB)$ is lc for some $t>0$. 
Then, there exists a Stein open subset $Y' \subset Y$ containing $W$ and a sequence of steps of a $(K_{X}+\Delta)$-MMP over $Y'$ around $W$ with scaling of a $\pi$-ample $\mathbb{R}$-divisor
$$(X,\Delta)=:(X_{0},\Delta_{0})\dashrightarrow (X_{1},\Delta_{1})\dashrightarrow \cdots \dashrightarrow (X_{n},\Delta_{n}),$$
which is represented by bimeromorphic contractions over $Y'$, such that $(X_{n},\Delta_{n})$ is a good minimal model or a Mori fiber space of $(X,\Delta)$ over $Y'$ around $W$. 
\end{thm} 

\begin{thm}[cf.~{\cite{hashizumehu}}]\label{thm--mmp-ampleeffective-intro} Let $\pi \colon X \to Y$ be a projective morphism from a normal analytic variety $X$ to a Stein space $Y$, and let $W \subset Y$ be a compact subset such that $\pi$ and $W$ satisfy (P). 
Let $(X,\Delta)$ be an lc pair and $A$ an effective $\pi$-ample $\mathbb{R}$-divisor on $X$ such that $(X,\Delta+A)$ is lc. 
Then, there exist a Stein open subset $Y' \subset Y$ containing $W$ and a sequence of steps of a $(K_{X}+\Delta+A)$-MMP over $Y'$ around $W$ with scaling of a $\pi$-ample $\mathbb{R}$-divisor
$$(X,\Delta+A)=:(X_{0},\Delta_{0}+A_{0})\dashrightarrow (X_{1},\Delta_{1}+A_{1})\dashrightarrow \cdots \dashrightarrow (X_{n},\Delta_{n}+A_{n}),$$
which is represented by bimeromorphic contractions over $Y'$, such that $(X_{n},\Delta_{n}+A_{n})$ is a good minimal model or a Mori fiber space of $(X,\Delta+A)$ over $Y'$ around $W$. 
\end{thm}

\begin{thm}\label{thm--strictmmp-2-intro} Let $\pi \colon X \to Y$ be a projective morphism from a normal analytic variety $X$ to a Stein space $Y$, and let $W \subset Y$ be a compact subset such that $\pi$ and $W$ satisfy (P). 
Let $(X,\Delta)$ be an lc pair such that $K_{X}+\Delta$ is pseudo-effective over $Y$. 
Let $A$ be a $\pi$-ample $\mathbb{R}$-divisor on $X$ such that $(X,\Delta+A)$ is lc and $K_{X}+\Delta+A$ is nef over $W$. 
Then, there exists a sequence of steps of a $(K_{X}+\Delta)$-MMP over $Y$ around $W$ with scaling of $A$
$$(X_{0},\Delta_{0})\dashrightarrow (X_{1},\Delta_{1})\dashrightarrow \cdots \dashrightarrow (X_{i},\Delta_{i}) \dashrightarrow \cdots$$
such that if we put
$$\lambda_{i}:={\rm inf}\{t \in \mathbb{R}_{\geq 0}\,|\, \text{$K_{X_{i}}+\Delta_{i}+tA_{i}$ {\rm is nef over} $W$}\}$$
for each $i \geq 0$, then ${\rm lim}_{i\to \infty}\lambda_{i}=0$. 
\end{thm}

By using Theorem \ref{thm--strictmmp-intro}, we prove a relation between the minimal model theory for projective lc pairs and that for complex analytic lc pairs. 

\begin{conj}\label{conj--algebraicmmp}
Let $(X,\Delta)$ be a projective lc pair. 
Then there exists a finite sequence of steps of a $(K_{X}+\Delta)$-MMP that terminates with a good minimal model or a Mori fiber space.   
\end{conj}

\begin{conj}\label{conj--analyticmmp}
Let $\pi \colon X \to Y$ be a projective morphism from a normal analytic variety $X$ to a Stein space $Y$, and let $W \subset Y$ be a compact subset such that $\pi$ and $W$ satisfy (P). 
Let $(X,\Delta)$ be an lc pair. 
Then, after shrinking $Y$around $W$, there exists a finite sequence of steps of a $(K_{X}+\Delta)$-MMP over $Y$ around $W$, which is represented by bimeromorphic contractions over $Y$, that terminates with a good minimal model or a Mori fiber space over $Y$ around $W$. 
\end{conj}

\begin{thm}\label{thm--mmp-equivalence-intro}
For every positive integer $d$, the following statements are equivalent:
\begin{enumerate}[(1)]
\item\label{thm--mmp-equivalence-intro-(1)}
Conjecture \ref{conj--algebraicmmp} holds for all projective klt pairs $(X,\Delta)$ such that ${\rm dim}\,X \leq d$. 
\item\label{thm--mmp-equivalence-intro-(2)}
Conjecture \ref{conj--algebraicmmp} holds for all projective lc pairs $(X,\Delta)$ such that ${\rm dim}\,X \leq d$. 
\item\label{thm--mmp-equivalence-intro-(3)}
Conjecture \ref{conj--analyticmmp} holds for all $\pi \colon X \to Y$, $W \subset Y$, and klt pairs $(X,\Delta)$ such that ${\rm dim}\,X \leq d$. 
\item\label{thm--mmp-equivalence-intro-(4)}
Conjecture \ref{conj--analyticmmp} holds for all $\pi \colon X \to Y$, $W \subset Y$, and lc pairs $(X,\Delta)$ such that ${\rm dim}\,X \leq d$. 
\end{enumerate}
\end{thm}

We also discuss MMP for algebraic stacks and analytic stacks. 
The definition of MMP for algebraic stacks and analytic stacks is different from that for algebraic varieties or complex analytic spaces (\cite[Section 2]{villalobos}, \cite[Definition 26.1]{lyu--murayama-mmp}, see also Definition \ref{outputMMP}). 
The main difficulty is gluing of MMP. 
This gluing problem was studied in \cite[Section 2]{villalobos} and \cite[Section 26]{lyu--murayama-mmp} for universally open and quasi-finite morphisms. 
In this paper, we consider the gluing problem for smooth morphisms.
By defining MMP for algebraic stacks and analytic stacks appropriately and using the gluing technique of the MMP, we obtain the following results.

\begin{thm}[Gluing version of Theorems~\ref{thm--strictmmp-intro}, \ref{thm--mmp-birkar-intro} and \ref{thm--mmp-ampleeffective-intro}]
\label{thm-glue}
Let $\pi \colon X \to Y$ be a projective morphism of algebraic stacks locally of finite type over $\C$ or complex analytic stacks.
Let $(X,\Delta)$ be an irreducible lc pair. 
Then the following holds:

\smallskip

\noindent
(1) A $(K_{X}+\Delta)$-MMP with scaling of a $\pi$-ample $\mathbb{R}$-line bundle $H$ over $Y$ exists.

\smallskip

\noindent
(2) Suppose one of the following conditions is satisfied:

\begin{enumerate}
\item[(a)]
$K_{X}+\Delta$ is not pseudo-effective over $Y$.

\item[(b)]
$(X,\Delta)$ has a log minimal model smooth locally on $Y$.
\end{enumerate}

\noindent
Then the MMP described in (1) terminates smooth locally on $Y$. 
The output of this MMP is a log minimal model or
a Mori fiber space of $(X,\Delta)$ over $Y$.

\smallskip

\noindent
(3) Suppose that $K_{X}+\Delta$ is pseudo-effective and one of the following conditions is satisfied:

\begin{enumerate}
   
\item [(c)]
$(X,\Delta)$ has a good minimal model smooth locally on $Y$.

\item[(d)]
Smooth locally on $Y$, there exists an effective $\mathbb{R}$-Cartier divisor $B$ on $X$
such that $K_{X}+\Delta+B\sim_{\mathbb{R},\,Y}0$ and $(X,\Delta+tB)$ is lc for some $t>0$. 

\item[(e)]
Smooth locally on $Y$, there exists an effective $\mathbb{R}$-Cartier divisor $A$ on $X$ which is ample over $Y$
such that $\Delta-A$ is effective.
\end{enumerate}

\noindent
Then the MMP described in (1) terminates smooth locally on $Y$. 
The output of this MMP is a good minimal model of $(X,\Delta)$ over $Y$.
\end{thm}

\begin{cor}\label{cor--local-criterion-mmp}
Let $\pi \colon X \to Y$ be a projective morphism of algebraic stacks locally of finite type over $\C$ or complex analytic stacks.
Let $(X,\Delta)$ be an irreducible  lc pair. 
Then the existence of a log minimal model (resp.\ a good minimal model, a Mori fiber space) of $(X,\Delta)$ over $Y$ can be checked smooth locally on $Y$.
\end{cor}

By taking $X$ and $Y$ in Theorem \ref{thm-glue} as complex analytic spaces, we obtain a result of the existence of a log minimal model or a Mori fiber space for complex analytic spaces without shrinking the base spaces (cf.~\cite[Section 27]{lyu--murayama-mmp}). 
We note that the definition of Mori fiber spaces for algebraic or analytic stacks (Definition \ref{defn--minmodel-stack}) is slightly different from that for algebraic varieties or complex analytic spaces (see Remark \ref{rem--MFS-difference}). 

The contents of this paper are as follows: 
In Section \ref{sec2}, we prove Theorems \ref{thm--strictmmp-intro}--\ref{thm--mmp-equivalence-intro}. 
In Section \ref{sec3}, we discuss MMP for algebraic stacks and analytic stacks, and we prove Theorem \ref{thm-glue} and Corollary \ref{cor--local-criterion-mmp}.

\begin{ack}
The first author was partially supported by JSPS KAKENHI Grant Number JP20K14297.
The second author was partially supported by JSPS KAKENHI Grant Number JP22K13887. 
The authors are grateful to Professor Osamu Fujino for discussions.   
We thank Professor Yoshinori Gongyo for a comment. 
We also thank the referee for carefully reading the manuscript and providing many valuable suggestions that greatly improved the paper. 
\end{ack}

\section{Minimal model program in complex analytic setting}\label{sec2}

Throughout this paper, {\em complex analytic spaces} are always assumed to be Hausdorff and second countable. 
{\em Analytic varieties} are reduced and irreducible complex analytic spaces. 
An $\mathbb{R}$-divisor on a normal complex variety is a (possibly infinite) formal sum of prime divisors that is locally finite (see, for example, \cite[Definition 2.32]{fujino-analytic-bchm}). 
We freely use notions of singularities of pairs in \cite[Section 3]{fujino-analytic-bchm}.

\begin{defn}[Property (P), see {\cite{fujino-analytic-bchm}}]\label{defn--property(P)}
Let $\pi \colon X \to Y$ be a projective morphism of complex analytic spaces, and let $W \subset Y$ be a compact subset. 
In this paper, we will use the following conditions:
\begin{itemize}
\item[(P1)]
$X$ is a normal analytic variety,
\item[(P2)]
$Y$ is a Stein space,
\item[(P3)]
$W$ is a Stein compact subset of $Y$, and
\item[(P4)]
$W \cap Z$ has only finitely many connected components for any analytic subset $Z$
which is defined over an open neighborhood of $W$. 
\end{itemize}
We say that {\em $\pi \colon X \to Y$ and $W \subset Y$ satisfy (P)} if the conditions (P1)--(P4) hold. 
\end{defn}

\begin{defn}[Models, {\cite[Definition 3.1]{eh-mmp}}]\label{defn--models}
Let $\pi \colon X \to Y$ be a projective morphism from a normal analytic variety $X$ to an analytic space $Y$, and let $(X,\Delta)$ be an lc pair. 
Let $W \subset Y$ be a subset. 
Let $\pi' \colon X' \to Y$ be a projective morphism from a normal analytic variety $X'$ to $Y$, and let $\phi \colon X \dashrightarrow X'$ be a bimeromorphic map over $Y$. 
We set $\Delta':=\phi_{*}\Delta+E$, where $E$ is the sum of all $\phi^{-1}$-exceptional prime divisors with all coefficients equal to $1$. 
When $K_{X'}+\Delta'$ is $\mathbb{R}$-Cartier, we say that $(X',\Delta')$ is a {\em log birational model of $(X,\Delta)$ over $Y$}. 

A log birational model of $(X,\Delta)$ over $Y$ is called a {\em weak log canonical model} ({\em weak lc model}, for short) {\em of $(X,\Delta)$ over $Y$ around $W$} if
\begin{itemize}
\item
$K_{X'}+\Delta'$ is nef over $W$, and
\item
for any prime divisor $P$ on $X$ that is exceptional over $X'$, we have 
$$a(P,X,\Delta) \leq a(P,X',\Delta').$$ 
\end{itemize} 
A weak lc model $(X',\Delta')$ of $(X,\Delta)$ over $Y$ around $W$ is a {\em log minimal model} if 
\begin{itemize}
\item
the above inequality on discrepancies is strict. 
\end{itemize}
A log minimal model $(X',\Delta')$ of $(X,\Delta)$ over $Y$ around $W$ is called a {\em good minimal model} if $K_{X'}+\Delta'$ is semi-ample over a neighborhood of $W$. 

Suppose that $W \subset Y$ is a compact subset such that $\pi \colon X \to Y$ and $W$ satisfy (P). 
Then, a log birational model $(X',\Delta')$ of $(X,\Delta)$ over $Y$ is called a {\em Mori fiber space over $Y$ around $W$} if there exists a contraction $X' \to Z$ over $Y$ such that
\begin{itemize}
\item
${\rm dim}\,X'>{\rm dim}\,Z$, 
\item
$\rho(X'/Y;W)-\rho(Z/Y;W)=1$ and $-(K_{X'}+\Delta')$ is ample over $Z$, and 
\item
for any prime divisor $P$ over $X$, we have 
$$a(P,X,\Delta) \leq a(P,X',\Delta'),$$
and the strict inequality holds if $P$ is a $\phi$-exceptional prime divisor on $X$.  
\end{itemize}

If $W=\emptyset$, then we formally define $(X,\Delta)$ itself to be a log minimal model (and thus a good minimal model) of $(X,\Delta)$ over $Y$ around $W$. 
Let $(\widetilde{X},\widetilde{\Delta})=\bigsqcup_{\lambda \in \Lambda}(\widetilde{X}_{\lambda},\widetilde{\Delta}_{\lambda})$ be a disjoint union of lc pairs, $\widetilde{\pi} \colon \widetilde{X} \to Y$ a projective morphism to an analytic space $Y$, and $W \subset Y$ a subset. 
Then a {\rm log birational model of $(\widetilde{X},\widetilde{\Delta})$ over $Y$} (resp.~a {\em log minimal model of $(\widetilde{X},\widetilde{\Delta})$ over $Y$ around $W$}) is a disjoint union of lc pairs $\bigsqcup_{\lambda \in \Lambda}(\widetilde{X}'_{\lambda},\widetilde{\Delta}'_{\lambda})$ such that $(\widetilde{X}'_{\lambda},\widetilde{\Delta}'_{\lambda})$ is a log birational model of $(\widetilde{X}_{\lambda},\widetilde{\Delta}_{\lambda})$ over $Y$ (resp.~a log minimal model of $(\widetilde{X}_{\lambda},\widetilde{\Delta}_{\lambda})$ over $Y$ around $W$) for all $\lambda \in \Lambda$. 
\end{defn}

\begin{defn}[Log MMP]\label{defn--mmp}
Let $\pi \colon X \to Y$ be a projective morphism from a normal analytic variety $X$ to a Stein space $Y$, and let $(X,\Delta)$ be an lc pair. 
Let $W \subset Y$ be a compact subset such that $\pi$ and $W$ satisfy (P).
\begin{enumerate}[(1)]
\item\label{defn--mmp-(1)}
A {\em step of a $(K_{X}+\Delta)$-MMP over $Y$ around $W$} is a diagram
 $$
\xymatrix{
(X,\Delta)\ar@{-->}[rr]^-{\phi}\ar[dr]^-{f}\ar[ddr]_-{\pi}&&(X',\Delta':=\phi_{*}\Delta)\ar[dl]_-{f'}\ar[ddl]^-{\pi'}\\
&Z\ar[d]\\
&Y
}
$$
consisting of projective morphisms such that 
\begin{itemize}
\item
$(X',\Delta')$ is an lc pair and $Z$ is a normal analytic variety,
\item
$\phi \colon X \dashrightarrow X'$ is a bimeromorphic contraction and $f\colon X \to Z$ and $f' \colon X' \to Z$ are bimeromorphic morphisms, 
\item
$f$ is a contraction of a $(K_{X}+\Delta)$-negative extremal ray of $\overline{\rm NE}(X/Y;W)$, in particular, $\rho(X/Y;W)-\rho(Z/Y;W)=1$ (cf.~\cite[Theorem 1.2 (4)(iii)]{fujino-analytic-conethm}) and $-(K_{X}+\Delta)$ is ample over $Z$, and
\item
$K_{X'}+\Delta'$ is ample over $Z$. 
\end{itemize}
Let $H$ be an $\mathbb{R}$-Cartier divisor on $X$. 
A {\em step of a $(K_{X}+\Delta)$-MMP over $Y$ around $W$ with scaling of $H$} is the above diagram satisfying the above conditions and
\begin{itemize}
\item
$K_{X}+\Delta+t H$ is nef over $W$ for some $t \in \mathbb{R}_{> 0}$, and 
\item
if we put 
$$\lambda:={\rm inf}\{t \in \mathbb{R}_{\geq 0}\,|\, \text{$K_{X}+\Delta+tH$ {\rm is nef over} $W$}\},$$
then $(K_{X}+\Delta+\lambda H)\cdot C= 0$ for any curve $C \subset \pi^{-1}(W)$ contracted by $f$. 
\end{itemize}

\item\label{defn--mmp-(2)}
A {\em sequence of steps of a $(K_{X}+\Delta)$-MMP over $Y$ around $W$} is a pair of sequences $\{Y_{i}\}_{i \geq 0}$ and $\{\phi_{i}\colon X_{i} \dashrightarrow X'_{i}\}_{i\geq 0}$, where $Y_{i} \subset Y$ are Stein open subsets and $\phi_{i}$ are bimeomorphic contractions of normal analytic varieties over $Y_{i}$, such that
\begin{itemize}
\item
$Y_{i}\supset W$ and $Y_{i}\supset Y_{i+1}$ for every $i\geq 0$, 
\item
$X_{0}=\pi^{-1}(Y_{0})$ and $X_{i+1}=X'_{i}\times_{Y_{i}}Y_{i+1}$ for every $i\geq 0$, and
\item
if we put $\Delta_{0}:=\Delta|_{X_{0}}$ and $\Delta_{i+1}:=(\phi_{i*}\Delta_{i})|_{X_{i+1}}$ for every $i\geq 0$, then 
$$(X_{i},\Delta_{i}) \dashrightarrow (X'_{i},\phi_{i*}\Delta_{i})$$ 
is a step of a $(K_{X_{i}}+\Delta_{i})$-MMP over $Y_{i}$ around $W$.  
\end{itemize}
For the simplicity of notation, a sequence of steps of a $(K_{X}+\Delta)$-MMP over $Y$ around $W$ is denoted by
$$(X_{0},\Delta_{0})\dashrightarrow (X_{1},\Delta_{1})\dashrightarrow \cdots \dashrightarrow (X_{i},\Delta_{i})\dashrightarrow \cdots.$$

\item\label{defn--mmp-(3)}
Let $H$ be an $\mathbb{R}$-Cartier divisor on $X$. 
We put 
$$H_{0}:=H|_{X_{0}} \qquad {\rm and} \qquad H_{i+1}:=(\phi_{i*}H_{i})|_{X_{i+1}}$$
for each $i \geq 0$. 
Then a {\em sequence of steps of a $(K_{X}+\Delta)$-MMP over $Y$ around $W$ with scaling of $H$} is a sequence of steps of a $(K_{X}+\Delta)$-MMP over $Y$ around $W$
$$(X_{0},\Delta_{0})\dashrightarrow (X_{1},\Delta_{1})\dashrightarrow \cdots \dashrightarrow (X_{i},\Delta_{i})\dashrightarrow \cdots$$
with the data $\{Y_{i}\}_{i \geq 0}$ and $\{\phi_{i}\colon X_{i} \dashrightarrow X'_{i}\}_{i\geq 0}$ such that $(X_{i},\Delta_{i})\dashrightarrow (X'_{i},\phi_{i*}\Delta_{i})$ is a step of a $(K_{X_{i}}+\Delta_{i})$-MMP over $Y$ around $W$ with scaling of $H_{i}$ for every $i$.
\end{enumerate}
\end{defn}

\begin{defn}
With notation as in Definition \ref{defn--mmp}, let
$$(X_{0},\Delta_{0})\dashrightarrow (X_{1},\Delta_{1})\dashrightarrow \cdots \dashrightarrow (X_{i},\Delta_{i})\dashrightarrow \cdots$$
be a sequence of steps of a $(K_{X}+\Delta)$-MMP over $Y$ around $W$ defined with $\{Y_{i}\}_{i \geq 0}$ and $\{\phi_{i}\colon X_{i} \dashrightarrow X'_{i}\}_{i\geq 0}$. 
We say that the $(K_{X}+\Delta)$-MMP is {\em represented by bimeromorphic contractions over $Y$} if $Y_{i}=Y$ for all $i \geq 0$. 
\end{defn}

We recall an important result proved by Fujino \cite{fujino-analytic-lcabundance}. 

\begin{thm}[cf.~{\cite[Theorem 1.7]{fujino-analytic-lcabundance}}]\label{thm--analytic-flip}
Let $\varphi \colon X \to Z$ be a projective bimeromorphic morphism of normal analytic varieties, and let $(X,\Delta)$ be an lc pair such that $-(K_{X}+\Delta)$ is $\varphi$-ample. 
Then we have a commutative diagram 
$$
\xymatrix{
(X,\Delta)
\ar[dr]_{\varphi}\ar@{-->}[rr]^{\phi} &&(X^{+},\Delta^{+})
\ar[dl]^{\varphi^{+}}\\
&Z
}
$$
satisfying the following properties: 
\begin{enumerate}[(i)]
\item
$\varphi^{+} \colon X^{+} \to Z$ is a projective small bimeromorphic morphism of normal analytic varieties, 
\item
$\Delta^{+}=\phi_{*}\Delta$ and $(X^{+},\Delta^{+})$ is lc, and
\item
$K_{X^{+}}+\Delta^{+}$ is $\varphi^{+}$-ample. 
\end{enumerate}
\end{thm}

\begin{rem}
In \cite[Theorem 1.7]{fujino-analytic-lcabundance}, the bimeromorphic morphism $\varphi \colon X \to Z$ is assumed to be small. 
However, this condition is not necessary in the proof. 
In fact, in the proof of \cite[Theorem 1.7]{fujino-analytic-lcabundance}, we first take a dlt blow-up, and then we construct $(X^{+},\Delta^{+})$ by using an MMP over $Z$ and the abundance. 
As in the algebraic case, it is easy to check that $(X^{+},\Delta^{+})$ satisfies the three properties.  
\end{rem}

The following lemma is the main technical result of this paper. 
This lemma plays an important role in the proofs of the main results. 

\begin{lem}[cf.~{\cite[Proposition 6.2]{hashizumehu}}]\label{lem--strict-decrease-nefthreshold} Let $\pi \colon X \to Y$ be a contraction from a normal analytic variety $X$ to a Stein space $Y$, and let $W \subset Y$ be a connected compact subset such that $\pi$ and $W$ satisfy (P). 
Let $(X,\Delta)$ be an lc pair. 
Suppose that $K_{X}+\Delta$ is not pseudo-effective over $Y$ or $(X,\Delta)$ has a log minimal model over $Y$ around $W$ after shrinking $Y$ around $W$. 
Then, there exist a Stein open subset $Y' \subset Y$ containing $W$ and a sequence of steps of a $(K_{X}+\Delta)$-MMP over $Y'$ around $W$ 
$$(X,\Delta)=:(X_{0},\Delta_{0})\dashrightarrow (X_{1},\Delta_{1})\dashrightarrow \cdots \dashrightarrow (X_{n},\Delta_{n}),$$
which is represented by bimeromorphic contractions over $Y'$, such that $(X_{n},\Delta_{n})$ is a log minimal model or a Mori fiber space of $(X,\Delta)$ over $Y'$ around $W$. 
\end{lem}

\begin{proof} 
The proof of \cite[Proposition 6.2]{hashizumehu} works in our situation. 

Suppose that some $(K_{X}+\Delta)$-MMP over $Y$ around $W$ contracts a prime divisor. In other words, suppose that there exists a finite sequence of steps of a $(K_{X}+\Delta)$-MMP over $Y$ around $W$ 
$$(X'_{0},\Delta'_{0})\dashrightarrow (X'_{1},\Delta'_{1})\dashrightarrow \cdots \dashrightarrow (X'_{m},\Delta'_{m})$$
such that after shrinking $Y$ around $W$, the bimeromorphic map $X'_{m-1}\dashrightarrow X'_{m}$ contracts a prime divisor $P_{m-1}$ on $X'_{m-1}$ whose image on $Y$ intersects $W$. 
By the lift of MMP in \cite[Subsection 3.7]{eh-mmp}, we get 
$$
\xymatrix{
(\widetilde{X}'_{0},\widetilde{\Delta}'_{0})\ar[d]_{g'_{0}}\ar@{-->}[r] &(\widetilde{X}'_{k_{1}},\widetilde{\Delta}'_{k_{1}})\ar[d]_{g'_{1}}\ar@{-->}[r]& \cdots\ar@{-->}[r]&  (\widetilde{X}'_{k_{m}},\widetilde{\Delta}'_{k_{m}})\ar[d]^{g'_{m}}\\
(X'_{0},\Delta'_{0})\ar@{-->}[r]&(X'_{1},\Delta'_{1})\ar@{-->}[r]& \cdots\ar@{-->}[r]&  (X'_{m},\Delta'_{m})
}
$$
such that 
\begin{itemize}
\item
the sequence of upper horizontal arrows is a sequence of steps of a $(K_{\widetilde{X}'_{0}}+\widetilde{\Delta}'_{0})$-MMP over $Y$ around $W$, and 
\item
each $g'_{i} \colon \widetilde{X}'_{k_{i}} \to X_{i}$ is a dlt blow-up of $(X'_{i},\Delta'_{i})$ and $\widetilde{X}'_{k_{i}}$ are $\mathbb{Q}$-factorial over $W$. 
\end{itemize}

As in the algebraic case, we have 
$$\rho(\widetilde{X}'_{k_{m}}/Y;W)<\rho(\widetilde{X}'_{k_{m-1}}/Y;W)\leq \rho(\widetilde{X}'_{0}/Y;W),$$ and we have $\rho(\widetilde{X}'_{0}/Y;W)<\infty$ by \cite[Subsection 4.1]{fujino-analytic-bchm}. 
We may replace $(X,\Delta)$ by $(X'_{m},\Delta'_{m})$ without loss of generality. 
After the replacement, if there exists a sequence of steps of a $(K_{X}+\Delta)$-MMP over $Y$ around $W$ that contracts a prime divisor whose image on $Y$ intersects $W$, then we can apply the above argument. 
Repeating this discussion, we may assume that any sequence of steps of a $(K_{X}+\Delta)$-MMP over $Y$ around $W$ does not contract any prime divisor whose image on $Y$ intersects $W$. 

By \cite[Theorem 3.9]{eh-mmp}, for any sequence of steps of a $(K_{X}+\Delta)$-MMP over $Y$ around $W$
 $$
(X''_{0},\Delta''_{0})\dashrightarrow (X''_{1},\Delta''_{1})\dashrightarrow \cdots \dashrightarrow (X''_{i},\Delta''_{i})\dashrightarrow \cdots,
$$
the strict transform of divisors defines a linear map
$$N^{1}(X''_{i-1}/Y;W) \longrightarrow N^{1}(X''_{i}/Y;W)$$
for each $i>0$, and the conclusion in the previous paragraph implies that the linear map is injective. 
Suppose that $\rho(X''_{i-1}/Y;W)<\rho(X''_{i}/Y;W)$ for some $i>0$. 
By using \cite[Subsection 3.7]{eh-mmp}, we get 
$$
\xymatrix{
(\widetilde{X}''_{0},\widetilde{\Delta}''_{0})\ar[d]_{g''_{0}}\ar@{-->}[r] &(\widetilde{X}''_{l_{1}},\widetilde{\Delta}''_{l_{1}})\ar[d]_{g''_{1}}\ar@{-->}[r]& \cdots\ar@{-->}[r]&  (\widetilde{X}''_{l_{i}},\widetilde{\Delta}''_{l_{i}})\ar[d]^{g''_{i}}\\
(X''_{0},\Delta''_{0})\ar@{-->}[r]&(X''_{1},\Delta''_{1})\ar@{-->}[r]& \cdots\ar@{-->}[r]&  (X''_{i},\Delta''_{i}),
}
$$
and we have 
$$0 \leq \rho(\widetilde{X}''_{l_{i}}/Y;W)-\rho(X''_{i}/Y;W) <  \rho(\widetilde{X}''_{0}/Y;W)-\rho(X''_{0}/Y;W)<\infty.$$
We may replace $(X,\Delta)$ by $(X''_{i},\Delta''_{i})$ without loss of generality, and we can repeat the above discussion. 
Repeating the discussion, we may assume that any sequence of steps of a $(K_{X}+\Delta)$-MMP over $Y$ around $W$
 $$
(X,\Delta)=:(X''_{0},\Delta''_{0})\dashrightarrow (X''_{1},\Delta''_{1})\dashrightarrow \cdots \dashrightarrow (X''_{i},\Delta''_{i})\dashrightarrow \cdots
$$
satisfies $\rho(X''_{i-1}/Y;W)=\rho(X''_{i}/Y;W)$ for all $i$. 
This implies that the linear map 
$$N^{1}(X''_{i-1}/Y;W) \longrightarrow N^{1}(X''_{i}/Y;W)$$
is bijective. 

By shrinking $Y$ around $W$, we may assume that we can write $$K_{X}+\Delta=\sum_{j=1}^{q}r_{j}D^{(j)}$$  for some $r_{j} \in \mathbb{R}$ and $\mathbb{Q}$-Cartier divisors $D^{(j)}$ on $X$. 
We put $\rho:=\rho(X/Y;W)$. We take $\alpha_{1},\cdots, \alpha_{\rho} \in \mathbb{R}_{>0}$ that are linearly independent over $\mathbb{Q}(r_{1},\cdots,\,r_{q})$, where $\mathbb{Q}(r_{1},\cdots,\,r_{q})$ is the field over $\mathbb{Q}$ generated by $r_{1},\cdots,\,r_{q}$. After shrinking $Y$ around $W$, we can find $\pi$-ample $\mathbb{Q}$-divisors $H^{(k)}$ ($1\leq k \leq \rho$) on $X$ such that \begin{itemize}\item $\{H^{(k)}\}^{\rho}_{k=1}$ is a basis of $N^{1}(X/Y;W)$, and \item $K_{X}+\Delta+\sum_{k=1}^{\rho}\alpha_{k}H^{(k)}$ is nef over $W$. \end{itemize} Set $H:=\sum_{k=1}^{\rho}\alpha_{k}H^{(k)}$. 
We run a $(K_{X}+\Delta)$-MMP over $Y$ around $W$ with scaling of $H$
$$(X_{0},\Delta_{0})\dashrightarrow (X_{1},\Delta_{1})\dashrightarrow \cdots \dashrightarrow (X_{j},\Delta_{j}) \dashrightarrow \cdots$$
and put
$$\lambda_{j}:={\rm inf}\{t \in \mathbb{R}_{\geq 0}\,|\, \text{$K_{X_{j}}+\Delta_{j}+tH_{j}$ {\rm is nef over} $W$}\}$$
for each $j \geq 0$. 
We show $\lambda_{j}>\lambda_{j+1}$ for all $j$. Suppose by contradiction that $\lambda_{i}=\lambda_{i+1}$ for some $i$. 
With notation as in \cite[Definition 3.5]{eh-mmp}, the $(i+1)$-th and $(i+2)$-th steps of the $(K_{X}+\Delta)$-MMP are written by $$ \xymatrix{ (X_{i},\Delta_{i})\ar@{-->}[rr]^{\phi_{i}}\ar[dr]_-{f_{i}}&&(X'_{i},\phi_{i*}\Delta_{i}),\ar[dl]^-{f'_{i}}\\ &Z_{i}\ar[d]\\ &{Y_{i}} }\quad \xymatrix{(X_{i+1},\Delta_{i+1})\ar@{-->}[rr]^{\phi_{i+1}}\ar[dr]_-{f_{i+1}}&&(X'_{i+1},\phi_{i+1*}\Delta_{i+1})\ar[dl]^-{f'_{i+1}}\\ &Z_{i+1}\ar[d]\\ &{Y_{i+1}} } $$ respectively such that $W \subset Y_{i+1}\subset Y_{i}\subset Y$ and $X_{i+1}=X'_{i}\times_{Y_{i}}Y_{i+1}$. By shrinking $Y$ around $W$, we may assume $Y_{i}=Y_{i+1}$. Then $X'_{i}=X_{i+1}$. Let $D^{(j)}_{i+1}$ (resp.~$H^{(k)}_{i+1}$, $H_{i+1}$) be a $\mathbb{Q}$-Cartier divisor on $X_{i+1}$ defined by the strict transforms of $D^{(j)}$ (resp.~$H^{(k)}$, $H$) and the restrictions repeatedly. By \cite[Theorem 3.9]{eh-mmp} and shrinking $Y$ around $W$, we may assume that all $D^{(j)}_{i+1}$ and $H^{(k)}_{i+1}$ are $\mathbb{Q}$-Cartier and $H_{i+1}$ is $\mathbb{R}$-Cartier.  Since $\rho(Z_{i}/Y;W)+1=\rho(X_{i}/Y;W)=\rho(X_{i+1}/Y;W)$, we see that $f'_{i}\colon X_{i+1}=X'_{i} \to Z_{i}$ is not a biholomorphism. Thus, there is a curve $C_{i}\subset X_{i+1}$ contracted by $f'_{i}$. Then $$(K_{X_{i+1}}+\Delta_{i+1}+\lambda_{i}H_{i+1})\cdot C_{i}=0\qquad {\rm and}\qquad (K_{X_{i+1}}+\Delta_{i+1})\cdot C_{i}>0.$$ We take a curve $C_{i+1} \subset X_{i+1}$ contracted by $f_{i+1}\colon X_{i+1}\to Z_{i+1}$. Then $$(K_{X_{i+1}}+\Delta_{i+1}+\lambda_{i+1}H_{i+1})\cdot C_{i+1}=0\qquad {\rm and}\qquad (K_{X_{i+1}}+\Delta_{i+1})\cdot C_{i+1}<0.$$ Then $$\frac{\lambda_{i+1}(H_{i+1}\cdot C_{i+1})}{\lambda_{i}(H_{i+1}\cdot C_{i})}=\frac{(K_{X_{i+1}}+\Delta_{i+1})\cdot C_{i+1}}{(K_{X_{i+1}}+\Delta_{i+1})\cdot C_{i}}.$$ Since $\lambda_{i}=\lambda_{i+1}$ and $K_{X_{i+1}}+\Delta_{i+1}=\sum_{j=1}^{q}r_{j}D_{i}^{(j)}$, putting $$\beta:=\frac{(K_{X_{i+1}}+\Delta_{i+1})\cdot C_{i+1}}{(K_{X_{i+1}}+\Delta_{i+1})\cdot C_{i}} \in \mathbb{Q}(r_{1},\cdots,\,r_{q}),$$ we have $$\sum_{k=1}^{\rho}\alpha_{k}H^{(k)}_{i+1}\cdot(C_{i+1}-\beta C_{i})=H_{i+1}\cdot (C_{i+1}-\beta C_{i})=0.$$ Since $H^{(k)}_{i+1}\cdot(C_{i+1}-\beta C_{i}) \in \mathbb{Q}(r_{1},\cdots,\,r_{q})$ for any $1 \leq k \leq \rho$ and $\alpha_{1},\cdots, \alpha_{\rho}$ are linearly independent over $\mathbb{Q}(r_{1},\cdots,\,r_{q})$, we have $H^{(k)}_{i+1}\cdot(C_{i+1}-\beta C_{i})=0$ for every $1 \leq k \leq \rho$. Moreover, since $\{H^{(k)}\}^{\rho}_{k=1}$ is a basis of $N^{1}(X/Y;W)$, it follows that $\{H^{(k)}_{i+1}\}^{\rho}_{k=1}$ generates $N^{1}(X_{i+1}/Y;W)$. Therefore, $C_{i+1}-\beta C_{i}=0$ as an element of $N_{1}(X/Y;W)$. This shows that $C_{i+1}$ and $C_{i}$ generate the same half line of $N_{1}(X/Y;W)$. But this is impossible because $(K_{X_{i+1}}+\Delta_{i+1})\cdot C_{i}>0$ and $(K_{X_{i+1}}+\Delta_{i+1})\cdot C_{i+1}<0$. Therefore, we have $\lambda_{i} \neq \lambda_{i+1}$. By \cite[Theorem 3.12]{eh-mmp}, we have $\lambda_{i}>\lambda_{i+1}$. 

Set $\lambda:={\rm lim}_{j\to \infty}\lambda_{j}$. 
By our assumption and \cite[Theorem 1.2]{eh-mmp}, $(X,\Delta+\lambda H)$ has a log minimial model over $Y$ around $W$. 
Then the $(K_{X}+\Delta)$-MMP over $Y$ around $W$
$$(X_{0},\Delta_{0})\dashrightarrow (X_{1},\Delta_{1})\dashrightarrow \cdots \dashrightarrow (X_{j},\Delta_{j}) \dashrightarrow \cdots$$
terminates by \cite[Theorem 3.15]{eh-mmp}. 
By shrinking $Y$ around $W$, we get a sequence of steps of a $(K_{X}+\Delta)$-MMP over $Y$ around $W$ 
$$(X_{0},\Delta_{0})\dashrightarrow (X_{1},\Delta_{1})\dashrightarrow \cdots \dashrightarrow (X_{n},\Delta_{n}),$$
which is represented by bimeromorphic contractions over $Y$, such that $(X_{n},\Delta_{n})$ is a log minimal model or a Mori fiber space of $(X,\Delta)$ over $Y$ around $W$. 
\end{proof}

We are ready to prove the main results of this paper. 
While Lemma \ref{lem--strict-decrease-nefthreshold} only shows the existence of a sequence of steps of an MMP that terminates with a log minimal model or a Mori fiber space over the base Stein space, the MMP in Theorem \ref{thm--strictmmp-intro}, Theorem \ref{thm--mmp-birkar-intro}, and Theorem \ref{thm--mmp-ampleeffective-intro} are MMP with scaling of relative ample $\mathbb{R}$-divisors.

\begin{proof}[Proof of Theorem \ref{thm--strictmmp-intro}]
We follow \cite[Proof of Theorem 1.7]{hashizumehu}. 
We use the notation as in Theorem \ref{thm--strictmmp-intro}. 
Suppose that $K_{X}+\Delta$ is not pseudo-effective over $Y$ or $(X,\Delta)$ has a log minimal model over $Y$ around $W$ after shrinking $Y$ around $W$. 

Put $(X_{0},\Delta_{0}):=(X,\Delta)$, $A_{0}:=A$, and 
$$\lambda_{0}:={\rm inf}\set{t \in \mathbb{R}_{\geq0}| \text{$K_{X_{0}}+\Delta_{0}+tA_{0}$ is nef over $W$}}.$$
If $\lambda_{0}=0$, then $(X_{0},\Delta_{0})$ is a log minimal model of $(X,\Delta)$ over $Y$ around $W$. 
If $\lambda_{0}>0$, then we pick $\lambda'_{0}\in (0,\lambda_{0})$ sufficiently close to $\lambda_{0}$ and we run a $(K_{X_{0}}+\Delta_{0}+\lambda'_{0}A_{0})$-MMP over $Y$ around $W$. 
Since $\lambda'_{0}>0$, \cite[Theorem 1.1]{eh-mmp} and Lemma \ref{lem--strict-decrease-nefthreshold} imply that this MMP terminates after finitely many steps. 
Therefore, we get a finite sequence of steps of a $(K_{X_{0}}+\Delta_{0}+\lambda'_{0}A_{0})$-MMP over $Y$ around $W$
$$(X_{0},\Delta_{0}+\lambda'_{0}A_{0}) \dashrightarrow (X_{1},\Delta_{1}+\lambda'_{0}A_{1}),$$
which is represented by a bimeromorphic contractions over some Stein open subset $Y_{1}$ of $Y$, such that $(X_{1},\Delta_{1}+\lambda'_{0}A_{1})$ is a log minimal model or a Mori fiber space over $Y_{1}$ around $W$.
By the length of extremal rays \cite[Theorem 1.1.6 (5)]{fujino-analytic-conethm} and choosing $\lambda'_{0}$ appropriately, we may assume that the strict transform of $K_{X_{0}}+\Delta_{0}+\lambda_{0}A_{0}$ trivially intersects the curve contracted in each step of the MMP. 
Thus, we may assume that this MMP is a $(K_{X_{0}}+\Delta_{0}+\lambda'_{0}A_{0})$-MMP with scaling of $A_{0}$ (Definition \ref{defn--mmp} (\ref{defn--mmp-(3)})). 
Then this MMP is a $(K_{X_{0}}+\Delta_{0})$-MMP with scaling of $A_{0}$ (cf.~\cite[Remark 3.7]{has-mmp-normal-pair}). 

If $(X_{1},\Delta_{1}+\lambda'_{0}A_{1})$ is a Mori fiber space of $(X_{0},\Delta_{0}+\lambda'_{0}A_{0})$ over $Y_{1}$ around $W$, then we complete the proof because $(X_{1},\Delta_{1})$ is a Mori fiber space of $(X_{0},\Delta_{0})$ over $Y_{1}$ around $W$. 
If $(X_{1},\Delta_{1}+\lambda'_{0}A_{1})$ is a log minimal model of $(X_{0},\Delta_{0}+\lambda'_{0}A_{0})$ over $Y_{1}$ around $W$, then we put 
$$\lambda_{1}:={\rm inf}\set{t \in \mathbb{R}_{\geq0}| \text{$K_{X_{1}}+\Delta_{1}+tA_{1}$ is nef over $W$}}.$$
Then $\lambda_{1}\lambda'_{0}<\lambda_{0}$. 
If $\lambda_{1}=0$, then $(X_{1},\Delta_{1})$ is a log minimal model of $(X,\Delta)$ over $Y$ around $W$. 
If $\lambda_{1}>0$, then we pick $\lambda'_{1}\in (0,\lambda_{1})$ sufficiently close to $\lambda_{1}$, and we apply the argument in the previous paragraph. 
Shrinking $Y$ around $W$, we get a finite sequence of steps of a $(K_{X_{1}}+\Delta_{1}+\lambda'_{1}A_{1})$-MMP over $Y_{1}$ around $W$ with scaling of $A_{1}$
$$(X_{1},\Delta_{1}+\lambda'_{1}A_{1}) \dashrightarrow (X_{2},\Delta_{2}+\lambda'_{1}A_{2}),$$
which is represented by a bimeromorphic contractions over some Stein open subset $Y_{2}$ of $Y_{1}$, such that $(X_{2},\Delta_{2}+\lambda'_{1}A_{2})$ is a log minimal model or a Mori fiber space over $Y_{2}$ around $W$.
By the length of extremal rays \cite[Theorem 1.1.6 (5)]{fujino-analytic-conethm} and choosing $\lambda'_{1}$ appropriately, we may assume that this MMP is a $(K_{X_{1}}+\Delta_{1}+\lambda'_{1}A_{1})$-MMP with scaling of $A_{1}$. 
Then 
$$(X_{0},\Delta_{0}) \dashrightarrow (X_{1},\Delta_{1}) \dashrightarrow (X_{2},\Delta_{2})$$
is a sequence of steps of a $(K_{X_{0}}+\Delta_{0})$-MMP over $Y$ around $W$ with scaling of $A_{0}$. 

By repeating the above discussion, we get a sequence of step of a $(K_{X}+\Delta)$-MMP over $Y$ around $W$ with scaling of $A$
$$(X,\Delta)=:(X_{0},\Delta_{0})\dashrightarrow (X_{1},\Delta_{1})\dashrightarrow \cdots \dashrightarrow (X_{i},\Delta_{i})\dashrightarrow \cdots$$
such that this MMP terminates with a Mori fiber space or $\lambda_{1}\neq {\rm lim}_{j \to \infty}\lambda_{j}$ for all $i$ when this MMP does not terminate. 
By \cite[Theorem 1.1]{eh-mmp} and  the assumption of Theorem \ref{thm--strictmmp-intro}, if $K_{X}+\Delta+\lambda A$ is pseudo-effective over $Y$ then $(X,\Delta+\lambda A)$ has a log minimal model over $Y$ around $W$. 
Thus, \cite[Theorem 3.15]{eh-mmp} shows that this MMP must terminate. 
From this argument, after shrinking $Y$ around $W$, we get a sequence of steps of a $(K_{X}+\Delta)$-MMP over $Y$ around $W$ with scaling of $A$
$$(X,\Delta)=:(X_{0},\Delta_{0})\dashrightarrow (X_{1},\Delta_{1})\dashrightarrow \cdots \dashrightarrow (X_{n},\Delta_{n}),$$
which is represented by bimeromorphic contractions over $Y$, such that $(X_{n},\Delta_{n})$ is a log minimal model or a Mori fiber space of $(X,\Delta)$ over $Y$ around $W$. 

Now suppose that $(X,\Delta)$ has a good minimal model $(X',\Delta')$ over $Y$ around $W$ after shrinking $Y$ around $W$. 
By shrinking $Y$ around $W$, we may assume that there is a common resolution $f_{n} \colon \tilde{X} \to X_{n}$ and $f' \colon \tilde{X} \to X'$ of the bimeromorphic map $X_{n} \dashrightarrow X'$ over $Y$. 
By \cite[Lemma 3.3]{eh-mmp}, we have 
$$f^{*}_{n}(K_{X_{n}}+\Delta_{n})=f'^{*}(K_{X'}+\Delta')$$
after shrinking $Y$ around $W$. 
Since $K_{X'}+\Delta'$ is semi-ample over a neighborhood of $W$, as in the algebraic case, we see that $K_{X_{n}}+\Delta_{n}$ is semi-ample over a neighborhood of $W$. 
Therefore, $(X_{n},\Delta_{n})$ is a good minimal model of $(X,\Delta)$ over $Y$ around $W$. 
We finish the proof. 
\end{proof}

\begin{proof}[Proof of Theorem \ref{thm--mmp-birkar-intro}]
By Theorem \ref{thm--strictmmp-intro}, it is sufficient to prove the existence of a good minimal model of $(X,\Delta)$ over $Y$ around $W$ when $K_{X}+\Delta$ is pseudo-effective over $Y$. 
By \cite[Lemma 2.16]{fujino-analytic-bchm}, we can find Stein open subsets $U_{1}$ and $U_{2}$ of $Y$ and a Stein compact subsets $W_{1}$ and $W_{2}$ of $Y$ 
such that 
$$W \subset U_{1} \subset W_{1} \subset U_{2} \subset W_{2}$$
and $\pi$ and $W_{l}$ satisfy (P) for $l=1,2$. 
By \cite[Theorem 1.2]{eh-mmp} and the construction of the MMP in \cite[Theorem 1.2]{eh-mmp}, after shrinking $Y$ around $W_{2}$, we get a log minimal model $(X',\Delta')$ of $(X,\Delta)$ over $Y$ around $W_{2}$ and an effective $\mathbb{R}$-Cartier divisors $B'$ on $X'$ such that $K_{X'}+\Delta'+B'\sim_{\mathbb{R},\,Y}0$ and $(X',\Delta'+t'B')$ is lc for some $t \gg t'>0$. 
By replacing $Y$ with $U_{2}$, we may assume that $K_{X'}+\Delta'$ is nef over $Y$. 
We also see that the numerical dimension of $(K_{X'}+\Delta')|_{S'^\nu}$ over $Y$ is zero or $-\infty$ for $S'^{\nu}:=X'$ or any lc center $S'$ of $(X',\Delta')$ with the normalization $S'^\nu \to S'$ (cf.~\cite[Remark 3.7]{has-finite}). 
By \cite{gongyo}, we see that $K_{X'}+\Delta'$ is log abundant over $Y$. 
By \cite[Theorem 1.5]{fujino-analytic-lcabundance}, we see that $K_{X'}+\Delta'$ is semi-ample over a neighborhood of $W$. 
This shows that $(X',\Delta')$ is a good minimal model of $(X,\Delta)$ over $Y$ around $W$. 
Thus, Theorem \ref{thm--mmp-birkar-intro} holds. 
\end{proof}

\begin{proof}[Proof of Theorem \ref{thm--mmp-ampleeffective-intro}]
By Theorem \ref{thm--strictmmp-intro}, it is sufficient to prove the existence of a good minimal model of $(X,\Delta+A)$ over $Y$ around $W$ when $K_{X}+\Delta+A$ is pseudo-effective over $Y$. 
By \cite[Lemma 2.16]{fujino-analytic-bchm}, we can find Stein open subsets $U_{1}$ and $U_{2}$ of $Y$ and a Stein compact subsets $W_{1}$ and $W_{2}$ of $Y$ 
such that 
$$W \subset U_{1} \subset W_{1} \subset U_{2} \subset W_{2}$$
and $\pi$ and $W_{l}$ satisfy (P) for $l=1,2$. 
By \cite[Theorem 1.1]{eh-mmp} and the construction of the MMP in \cite[Theorem 1.1]{eh-mmp}, after shrinking $Y$ around $W_{2}$, we get a log minimal model $(X',\Gamma')$ of $(X,\Delta+A)$ over $Y$ around $W_{2}$.  
By replacing $Y$ with $U_{2}$, we may assume that $K_{X'}+\Gamma'$ is nef over $Y$. 
For $S'^{\nu}:=X'$ or any lc center $S'$ of $(X',\Gamma')$ with the normalization $S'^{\nu} \to S'$, the restriction of $(K_{X'}+\Gamma')|_{S'^{\nu}}$ to the general fiber of the Stein factorization of $S'^{\nu} \to Y'$ is semi-ample (\cite[Theorem 1.5]{hashizumehu}). 
Thus, $K_{X'}+\Gamma'$ is log abundant over $Y$ (see also \cite[Remark 3.7]{has-finite}). 
Then \cite[Theorem 1.5]{fujino-analytic-lcabundance} implies that $K_{X'}+\Gamma'$ is semi-ample over a neighborhood of $W$. 
Therefore, $(X',\Gamma')$ is a good minimal model of $(X,\Delta)$ over $Y$ around $W$. 
Thus, Theorem \ref{thm--mmp-ampleeffective-intro} holds. 
\end{proof}

\begin{proof}[Proof of Theorem \ref{thm--strictmmp-2-intro}]
As in \cite[Proof of Theorem 1.7]{hashizumehu} or \cite[Corollary 3.9]{has-finite}, we can construct a sequence of steps of a $(K_{X}+\Delta)$-MMP over $Y$ around $W$ with scaling of $A$
$$(X,\Delta)=:(X_{0},\Delta_{0})\dashrightarrow (X_{1},\Delta_{1})\dashrightarrow \cdots \dashrightarrow (X_{i},\Delta_{i}) \dashrightarrow \cdots$$
such that if we put
$$\lambda_{i}:={\rm inf}\{t \in \mathbb{R}_{\geq 0}\,|\, \text{$K_{X_{i}}+\Delta_{i}+tA_{i}$ {\rm is nef over} $W$}\}$$
for each $i \geq 0$ and $\lambda:={\rm lim}_{i\to \infty}\lambda_{i}$, then $\lambda \neq \lambda_{j}$ for all $j$. 
It is enough to prove $\lambda =0$. 
Suppose by contradiction that $\lambda > 0$. 
Then the above MMP is also a sequence of steps of a $(K_{X}+\Delta+\lambda A)$-MMP over $Y$ around $W$ with scaling of $(1-\lambda) A$
$$(X,\Delta+\lambda A)=:(X_{0},\Delta_{0}+\lambda A_{0})\dashrightarrow  \cdots \dashrightarrow (X_{i},\Delta_{i}+\lambda A_{i}) \dashrightarrow \cdots$$
such that if we put
$$\mu_{i}:={\rm inf}\{t \in \mathbb{R}_{\geq 0}\,|\, \text{$K_{X_{i}}+\Delta_{i}+(t+\lambda)A_{i}$ {\rm is nef over} $W$}\}$$
for each $i \geq 0$, then ${\rm lim}_{i\to \infty}\mu_{i}=0$. 
By \cite[Theorem 3.15]{eh-mmp}, we have $\mu_{m}=0$ for some $m \geq 0$, equivalently, $\lambda = \lambda_{m}$ for some $m \geq 0$, contradicting the fact that $\lambda \neq \lambda_{j}$ for all $j$. 
Thus, we have $\lambda =0$, and we complete the proof. 
\end{proof}

\begin{proof}[Proof of Theorem \ref{thm--mmp-equivalence-intro}]
It is obvious that (\ref{thm--mmp-equivalence-intro-(4)}) implies (\ref{thm--mmp-equivalence-intro-(3)}), and (\ref{thm--mmp-equivalence-intro-(3)}) implies (\ref{thm--mmp-equivalence-intro-(1)}) because (\ref{thm--mmp-equivalence-intro-(1)}) is a special case of (\ref{thm--mmp-equivalence-intro-(3)}) where the base space is a point. 
It is known by \cite{fg-lcring} that (\ref{thm--mmp-equivalence-intro-(1)}) implies (\ref{thm--mmp-equivalence-intro-(2)}). 
Therefore, it is enough to show that (\ref{thm--mmp-equivalence-intro-(2)}) implies (\ref{thm--mmp-equivalence-intro-(4)}). 

Assume the condition (\ref{thm--mmp-equivalence-intro-(2)}), in other words, assume  Conjecture \ref{conj--algebraicmmp} for all projective lc pairs of dimension $\leq d$. 
Let $\pi \colon X \to Y$, $W \subset Y$, and $(X,\Delta)$ be as in Conjecture \ref{conj--analyticmmp} such that ${\rm dim}\,X \leq d$. 
If $K_{X}+\Delta$ is not $\pi$-pseudo-effective, then Theorem \ref{thm--strictmmp-intro} implies that (\ref{thm--mmp-equivalence-intro-(4)}) holds for $\pi \colon X \to Y$, $W \subset Y$, and $(X,\Delta)$. 
Therefore, we may assume that $K_{X}+\Delta$ is $\pi$-pseudo-effective. 
By Theorem \ref{thm--strictmmp-intro} again, it is enough to prove the existence of a good minimal model of $(X,\Delta)$ over $Y$ around $W$. 
By \cite[Lemma 2.16]{fujino-analytic-bchm}, we can find Stein open subsets $U_{1}$ and $U_{2}$ of $Y$ and a Stein compact subsets $W_{1}$ and $W_{2}$ of $Y$ 
such that 
$$W \subset U_{1} \subset W_{1} \subset U_{2} \subset W_{2}$$
and $\pi$ and $W_{l}$ satisfy (P) for $l=1,2$. 
By Theorem \ref{thm--strictmmp-2-intro}, there exists a sequence of steps of a $(K_{X}+\Delta)$-MMP over $Y$ around $W_{2}$ with scaling of a $\pi$-ample $\mathbb{R}$-divisor $A$
$$(X,\Delta)=:(X_{0},\Delta_{0})\dashrightarrow (X_{1},\Delta_{1})\dashrightarrow \cdots \dashrightarrow (X_{i},\Delta_{i}) \dashrightarrow \cdots$$
such that if we put
$$\lambda_{i}:={\rm inf}\{t \in \mathbb{R}_{\geq 0}\,|\, \text{$K_{X_{i}}+\Delta_{i}+tA_{i}$ {\rm is nef over} $W_{2}$}\}$$
for each $i \geq 0$, then ${\rm lim}_{i\to \infty}\lambda_{i}=0$. 
Since we assume Conjecture \ref{conj--algebraicmmp} for all projective lc pairs of dimension $\leq d$, we see that all the lc pairs $(X_{i},\Delta_{i})$ are log abundant over $Y$ around $W_{2}$. 
Indeed, over a neighborhood of $W_{2}$ we take a dlt blow-up $(\tilde{X}_{i},\tilde{\Delta}_{i})$ of $(X_{i},\Delta_{i})$. 
For $\tilde{S}_{i}:=\tilde{X}_{i}$ or any lc center $\tilde{S}_{i}$ of $(\tilde{X}_{i},\tilde{\Delta}_{i})$, let $(\tilde{S}_{i}, \Delta_{\tilde{S}_{i}})$ be the lc pair defined by adjunction $K_{\tilde{S}_{i}}+\Delta_{\tilde{S}_{i}}=(K_{\tilde{X}_{i}}+\tilde{\Delta}_{i})|_{\Delta_{\tilde{S}_{i}}}$. 
By Conjecture \ref{conj--algebraicmmp} for projective lc pairs of dimension $\leq d$, the restriction of $(\tilde{S}_{i}, \Delta_{\tilde{S}_{i}})$ to the general fiber of the Stein factorization of $\tilde{S}_{i} \to Y$ has a good minimal model or a Mori fiber space. 
Thus, $K_{\tilde{S}_{i}}+\Delta_{\tilde{S}_{i}}$ is abundant over $Y$ around $W_{2}$, and therefore $(X_{i},\Delta_{i})$ is log abundant over $Y$ around $W_{2}$ for all $i$. 
By \cite[Theorem 1.3]{eh-mmp}, this MMP terminates after finitely many steps. 
Let $(X_{n},\Delta_{n})$ be the resulting lc pair. 
Since $K_{X}+\Delta$ is $\pi$-pseudo-effective, $K_{X_{n}}+\Delta_{n}$ is nef over $W_{2}$. 
By replacing $Y$ with $U_{2}$, we may assume that $K_{X_{n}}+\Delta_{n}$ is nef over $Y$. 
Now 
$$W \subset U_{1} \subset W_{1} \subset Y$$
and the abundance conjecture holds for all projective lc pairs of dimension $\leq d$ becasue we assume that Conjecture \ref{conj--algebraicmmp} holds for all projective lc pairs of dimension $\leq d$. 
By \cite[Corollary 1.19]{fujino-analytic-lcabundance}, we see that $K_{X_{n}}+\Delta_{n}$ is semi-ample over a neighborhood of $W$. 
Hence, $(X_{n},\Delta_{n})$ is a good minimal model over $Y$ around $W$, and therefore (\ref{thm--mmp-equivalence-intro-(4)}) holds by Theorem \ref{thm--strictmmp-intro}. 

By the above argument, we see that (\ref{thm--mmp-equivalence-intro-(2)}) implies (\ref{thm--mmp-equivalence-intro-(4)}). 
Thus, Theorem \ref{thm--mmp-equivalence-intro} holds. 
\end{proof}

\section{Gluing of MMP}\label{sec3}
In this section, we discuss MMP for algebraic stacks and analytic stacks. 
We refer to \cite{villalobos} and \cite{lyu--murayama-mmp}. 
In this section we usually use globally $\mathbb{R}$-Cartier divisors on complex analytic varieties, though we implicitly used the notion in the previous section. 
Recall that a {\em globally $\mathbb{R}$-Cartier divisor} on a complex analytic variety is a finite $\mathbb{R}$-linear combination of Cartier divisors (\cite[Definition 2.32]{fujino-analytic-bchm}). 

We first construct a variant of sequences of steps of MMP. 
We follow the construction in \cite{villalobos} and \cite[Definition 1]{kollar-non-Qfac-mmp}. 
We only discuss the complex analytic case because the algebraic case is similar and simpler. 

\begin{defn}\label{defn--mmp-another}
Let $\pi \colon X \to Y$ be a projective morphism from a normal analytic variety $X$ to a Stein space $Y$, and let $W \subset Y$ be a compact subset such that $\pi$ and $W$ satisfy (P). 
Let $(X,\Delta)$ be an lc pair. 
By shrinking $Y$ around $W$, we may assume that $K_{X}+\Delta$ is globally $\mathbb{R}$-Cartier.  
Let $H$ be an $\mathbb{R}$-Cartier divisor on $X$ such that $K_{X}+\Delta+c_{0}H$ is ample over $Y$ around $W$ for some $c_{0}>0$. 
We define 
$$\lambda_{0}:={\rm inf}\{\nu \in \mathbb{R}_{\geq 0}\,|\, \text{$K_{X}+\Delta+\nu H$ {\rm is nef over} $W$}\}.$$
Then $c_{0}>\lambda_{0}$ and $\frac{c_{0}}{c_{0}-\lambda_{0}}(K_{X}+\Delta+\lambda_{0}H)=K_{X}+\Delta+\frac{\lambda_{0}}{c_{0}-\lambda_{0}}(K_{X}+\Delta+c_{0}H)$. 
If $\lambda_{0}=0$, we stop the discussion. 
If $\lambda_{0}>0$, by the basepoint-free theorem (\cite[Theorem 5.3.1]{fujino-analytic-conethm}), there exist a Stein open neighborhood $Y_{0} \supset W$ and the contraction $\varphi \colon \pi^{-1}(Y_{0}) \to V_{0}$ over $Y_{0}$ defined by $K_{X}+\Delta+\lambda_{0}H$. 
We set $X_{0}:=\pi^{-1}(Y_{0})$, $\Delta_{0}:=\Delta|_{X_{0}}$, and $H_{0}:=H|_{X_{0}}$. 
By shrinking $Y_{0}$ around $W$, we may assume that $-(K_{X_{0}}+\Delta_{0})$ is ample over $Y_{0}$. 
If $K_{X_{0}}+\Delta_{0}+\lambda_{0}H_{0}$ is not big over $Y$, in other words, if ${\rm dim}\,X>{\rm dim}\,V_{0}$, then we stop the discussion. 
If $K_{X_{0}}+\Delta_{0}+\lambda_{0}H_{0}$ is big over $Y$, then Theorem \ref{thm--analytic-flip} implies that after shrinking $Y_{0}$ around $W$ suitably, we get a diagram
$$
\xymatrix{
(X_{0},\Delta_{0})
\ar[dr]_{\varphi}\ar@{-->}[rr]^{\phi} &&(X'_{0},\Delta'_{0})
\ar[dl]^{\varphi'}\\
&V_{0}
}
$$
over $Y_{0}$ satisfying the following properties: 
\begin{itemize}
\item
$\varphi' \colon X'_{0} \to V_{0}$ is a small projective bimeromorphic morphism of normal analytic varieties, 
\item
$\Delta'_{0}=\phi_{*}\Delta'_{0}$ and $(X'_{0},\Delta'_{0})$ is lc, and
\item
$K_{X'_{0}}+\Delta'_{0}$ is $\varphi'$-ample. 
\end{itemize}
Put $H'_{0}:=\phi_{*} H_{0}$. 
Since $\varphi \colon X_{0} \to V_{0}$ is the contraction over $Y_{0}$ defined by $K_{X}+\Delta+\lambda_{0}H$, the $\varphi'$-ampleness of $K_{X'_{0}}+\Delta'_{0}$ implies that $K_{X'_{0}}+\Delta'_{0}+(\lambda_{0}-\varepsilon)H'_{0}$ is ample over a neighborhood of $W$ for any $0< \varepsilon \ll 1$. 
By shrinking $Y_{0}$ around $W$ suitably, we can find $c_{1}<\lambda_{0}$ such that $K_{X'_{0}}+\Delta'_{0}+c_{1}H'_{0}$ is ample over $Y_{0}$. 

We  apply the above argument to $(X'_{0},\Delta'_{0}) \to Y_{0}$, and $H'_{0}$, and we may repeat the above argument as far as $K_{X_{i}}+\Delta_{i}+\lambda_{i}H_{i}$ is big over $Y$ in the $(i+1)$-th step of the argument. 
\end{defn}

The following result immediately follows from the above construction. 

\begin{thm}\label{thm-mmp-another}
Let $\pi \colon X \to Y$ be a projective morphism from a normal analytic variety $X$ to a Stein space $Y$, and let $W \subset Y$ be a compact subset such that $\pi$ and $W$ satisfy (P). 
Let $(X,\Delta)$ be an lc pair such that $K_{X}+\Delta$ is globally $\mathbb{R}$-Cartier.  
Let $H$ be a globally $\mathbb{R}$-Cartier divisor on $X$ such that $K_{X}+\Delta+c_{0}H$ is ample over a neighborhood of $W$. 
Then there exists a triple of sequences $\{Y_{i}\}_{i \geq 0}$, $\{\phi_{i}\colon X_{i} \dashrightarrow X'_{i}\}_{i\geq 0}$, and $\{\lambda_{i}\}_{i \geq 0}$, where $Y_{i} \subset Y$ are Stein open subsets and $\phi_{i}$ are bimeomorphic contractions of normal analytic varieties over $Y_{i}$, such that
\begin{itemize}
\item
$Y_{i}\supset W$ and $Y_{i}\supset Y_{i+1}$ for every $i\geq 0$, 
\item
$X_{0}=\pi^{-1}(Y_{0})$ and $X_{i+1}=X'_{i}\times_{Y_{i}}Y_{i+1}$ for every $i\geq 0$, and
\item
if we put $\Delta_{0}:=\Delta|_{X_{0}}$ and $\Delta_{i+1}:=(\phi_{i*}\Delta_{i})|_{X_{i+1}}$ for every $i\geq 0$, then 
\begin{itemize}
\item
$\lambda_{i}$ is defined to be
$$\lambda_{i}:={\rm inf}\{\nu \in \mathbb{R}_{\geq 0}\,|\, \text{$K_{X_{i}}+\Delta_{i}+\nu H_{i}$ {\rm is nef over} $W$}\}$$
and $\lambda_{i}>\lambda_{i+1}$ for each $i \geq 0$, 
\item
for each $i \geq 1$, the bimeromorphic contraction 
$\pi^{-1}(Y_{i}) \dashrightarrow X_{i}$ over $Y_{i}$ is the ample model of $(K_{X}+\Delta+(\lambda_{i-1}-\varepsilon)H)|_{\pi^{-1}(Y_{i})}$ over $Y_{i}$ {\rm (}Definition \ref{defn--stack-amplemodel} below,  
 {\rm cf.~{\cite[Definition~3.6.5]{bchm}}}{\rm )} for any $0<\varepsilon \ll 1$, and
\item
if $(\pi^{-1}(Y_{i}), (\Delta+\lambda_{i}H)|_{\pi^{-1}(Y_{i})})$ is lc, then the bimeromorphic contraction 
$$(\pi^{-1}(Y_{i}), (\Delta+\lambda_{i}H)|_{\pi^{-1}(Y_{i})})\dashrightarrow (X_{i},\Delta_{i}+\lambda_{i} H_{i})$$
 is a log minimal model over $Y$ around $W$. 
\end{itemize}
\end{itemize}
\end{thm}

\begin{rem}
The MMP in Definition \ref{defn--mmp-another} or Theorem \ref{thm-mmp-another} is also a $(K_{X}+\Delta)$-MMP with scaling of $K_{X}+\Delta+c_{0}H$ (see \cite[Remark 3.7]{has-mmp-normal-pair} for the algebraic case). 
Hence, we can often replace $H$ by $K_{X}+\Delta+c_{0}H$ and we may assume that $H$ is nef over $W$. 
\end{rem}

\begin{lem}\label{lem--mmp-another-termination}
With notation as in Theorem \ref{thm-mmp-another}, let $\mu$ be the pseudo-effective threshold of $K_{X}+\Delta$ with respect to $H$ over $Y$. 
Then $\lambda_{i} \geq \mu$ for all $i$ and ${\rm lim}_{i \to \infty}\lambda_{i}=\mu$.
Furthermore, if $(X,\Delta)$ has a log minimal model (resp.\ a good minimal model, a Mori fiber space) over $Y$ around $W$ after shrinking $Y$ around $W$, the MMP terminates after finitely many steps with a log minimal model (resp.\ a good minimal model, a Mori fiber space) over $Y$ around $W$. 
\end{lem}

\begin{proof}
Since $c_{0}>\lambda_{0}$ and $\frac{c_{0}}{c_{0}-t}(K_{X}+\Delta+tH)=K_{X}+\Delta+\frac{t}{c_{0}-t}(K_{X}+\Delta+c_{0}H)$ for any $0 \leq t <c_{0}$, we may replace $H$, $\mu$, and $\{\lambda_{i}\}_{i \geq 0}$ by $K_{X}+\Delta+c_{0}H$, $\frac{\mu}{c_{0}-\mu}$, and $\{\frac{\lambda_{i}}{c_{0}-\lambda_{i}}\}_{i \geq 0}$, respectively. 
Thus, we may assume that $H$ is ample over $Y$. 
In particular, for any $t>0$, after shrinking $Y$ around $W$, we may assume that $(X,\Delta+tH)$ is lc and the pair has a log minimal model or a Mori fiber space over $Y$ around $W$. 

We put $\lambda_{\infty}={\rm lim}_{i \to \infty}\lambda_{i}$. 
Assume to the contrary that $\lambda_{\infty} > \mu$. 
Then $(X,\Delta+\lambda_{\infty}H)$ has a good minimal model (Theorem \ref{thm--mmp-birkar-intro}), and we can get a contradiction if $\lambda_{m}=\lambda_{\infty}$ for some $m$ since $\lambda_{m}>\lambda_{m+1}$. 
Hence, to prove $\lambda_{\infty}=\mu$, we may replace $\Delta$ by $\Delta+\lambda_{\infty}H$.
Note that if the MMP terminates after finitely many steps, then the claim $\lambda_{\infty}=\mu$ clearly holds.
Therefore, to prove Theorem \ref{lem--mmp-another-termination}, it is sufficient to prove $\lambda_{m}=0$ for some $m$ when ${\rm lim}_{i \to \infty}\lambda_{i}=0$ and $(X,\Delta)$ has a log minimal model over $Y$ around $W$ after shrinking $Y$ around $W$. 

By Theorem \ref{thm--strictmmp-intro} and shrinking $Y$ around $W$, we have a bimeromorphic contraction $(X,\Delta)\dashrightarrow (X',\Delta')$ over $Y$ and a positive real number $c$ such that $K_{X'}+\Delta'+tH'$ is nef over $W$ for any $t \in [0,c]$. 
We pick an index $m$ such that $\lambda_{m}<c$. 
We show that $\lambda_{m}=0$. 
By shrinking $Y$, we may assume $Y=Y_{m}$, and we have a bimeromorphic contraction $X \dashrightarrow X_{m}$ over $Y$. 
Let $g \colon \tilde{X} \to X$ be a resolution such that the induced bimeromorphic maps $g_{m}\colon \tilde{X} \dashrightarrow X_{m}$ and $g' \colon \tilde{X} \dashrightarrow X'$ are morphisms. 
Since $\lambda_{m}<c$, there are at least two points, denoted by $t_{1}$ and $t_{2}$, in $[0,c]\cap [\lambda_{m},\lambda_{m-1}]$. 
By the negativity lemma (see \cite[Corollary 2.16]{eh-mmp}) and shrinking $Y$ around $W$, we have
\begin{equation*}
\begin{split}
{g_{m}}^{*}(K_{X_{m}}+\Delta_{m}+t_{1}H_{m})&=g'^{*}(K_{X'}+\Delta'+t_{1}H'), \quad {\rm and} \\
{g_{m}}^{*}(K_{X_{m}}+\Delta_{m}+t_{2}H_{m})&=g'^{*}(K_{X'}+\Delta'+t_{2}H').
\end{split}
\end{equation*}
Then ${g_{m}}^{*}(K_{X_{m}}+\Delta_{m})=g'^{*}(K_{X'}+\Delta')$. 
From this we see that $K_{X_{m}}+\Delta_{m}$ is nef over $W$, and therefore $\lambda_{m}=0$. 
Note that $K_{X_{m}}+\Delta_{m}$ is semi-ample over a neighborhood of $W$ if and only if so is $K_{X'}+\Delta'$.
Hence Theorem \ref{lem--mmp-another-termination} holds. 
\end{proof}

From now on, we discuss the gluing problem of the above sequences to define MMP and resulting models for algebraic stacks and analytic stacks.  
All algebraic stacks are assumed to be locally of finite type over a base field $k$ (More generally, the gluing results in this section work for excellent algebraic stacks $X$ admitting a dualizing complex $\omega_{X}^{\bullet}$. 
See Remark~\ref{rem:logcanodiv}).
{\em Complex analytic stacks} are defined as stacks in groupoids $X$ over the category of complex analytic spaces equipped with the \'{e}tale topology such that
the diagonal map $\Delta_{X}$ is representable by complex analytic spaces and there exists a smooth covering $U\to X$ from a complex analytic space $U$ (for the definition of stacks in groupoids, see {\cite[02ZI]{Stacks}}).
All arguments in this section are conducted within the category of schemes or complex analytic spaces. The term {\em stacks} refers to either algebraic stacks or complex analytic stacks,
and {\em spaces} refers to algebraic spaces or complex analytic spaces unless otherwise stated.

\begin{defn}[Line bundle]
Let $X$ be a stack.
In this paper, a {\em line bundle} $\mathcal{L}$ on $X$ means an invertible $\O_{X}$-module on the lisse-\'{e}tale site of $X$.
Note that for complex analytic stacks, lisse-\'{e}tale site can also be defined as in algebraic case {\cite[0787]{Stacks}}. 
It can be regarded as a family of line bundles $\{\mathcal{L}_{U}\}_{U\to X}$ on spaces $U$ equipped with a smooth morphism $U\to X$ such that
for any morphism $\varphi\colon U\to V$ over $X$, 
an isomorphism $\rho_{\varphi}\colon \varphi^{*}\mathcal{L}_{V}\to \mathcal{L}_{U}$ is attached,
which is compatible with compositions of morphisms. 

The group of isomorphism classes of line bundles on $X$, which we call the Picard group of $X$, is identified with
$\mathrm{Pic}(X)=H^{1}(X,\O^{*}_{X})$.
For $R=\Q$ or $\R$, {\em an $R$-line bundle} on $X$ can also be defined by an element of 
$\mathrm{Pic}_{R}(X):=H^{1}(X,\O^{*}_{X}\otimes_{\Z}R)$.
Note that an $R$-line bundle may not come from an element of $\mathrm{Pic}(X)\otimes_{\Z}R$ while it comes smooth locally on $X$.
\end{defn}

\begin{defn}[Relative line bundle]
Let $\pi \colon X\to Y$ be a proper morphism of stacks and
let $R=\Z, \Q$ or $\R$.
The {\em $R$-Picard sheaf on $Y$}, which is denoted by $\mathrm{Pic}_{\pi, R}$, is defined as the sheafification of the presheaf in the smooth topology on $Y$
which sends a space $U$ smooth over $Y$ to $\mathrm{Pic}_{R}(X\times_{Y}U)/\pi^{*}\mathrm{Pic}_{R}(U)$.
A global section $\mathcal{L}$ of $\mathrm{Pic}_{\pi, R}$ is called a {\em (relative) $R$-line bundle on $X$ over $Y$}.
It has a representation 
$$
\mathcal{L}=(\{Y_{i}\to Y\}_{i}, \{\mathcal{L}_{i}\}_{i})
$$
for some smooth covering $\{Y_{i}\to Y\}_{i}$ from spaces $Y_i$ and some $R$-line bundles $\mathcal{L}_{i}$ on $X_i$ such that $\mathcal{L}_{i}|_{X_{ij}}$ and $\mathcal{L}_{j}|_{X_{ij}}$ are equal in $\mathrm{Pic}_{R}(X_{ij})/\pi^{*}\mathrm{Pic}_{R}(Y_{ij})$,
where $X_i:=X\times_{Y}Y_{i}$, $Y_{ij}:=Y_{i}\times_{Y}Y_{j}$ and $X_{ij}:=X\times_{Y}Y_{ij}$. 

For a relative $R$-line bundle $\mathcal{L}$ on $X$ over $Y$ and a morphism $Y' \to Y$ from a space $Y'$, the pullback of $\mathcal{L}$ to $X \times_{Y}Y'$ is denoted by $\mathcal{L}|_{X \times_{Y}Y'}$ or $\mathcal{L}|_{Y'}$ if there is no risk of confusion. 
\end{defn}

\begin{defn}[Ample line bundle]
Let $\pi \colon X\to Y$ be a proper morphism of stacks.
A line bundle $\mathcal{L}$ on $X$ is {\em basepoint-free} over $Y$ if the natural map $\pi^{*}\pi_{*}\mathcal{L}\to \mathcal{L}$ is surjective.
Then this surjection induces a morphism $\varphi_{\mathcal{L}}\colon X\to \mathbb{P}_{Y}(\pi_{*}\mathcal{L})$ by applying $\mathbb{P}_{Y}(-):=\mathrm{Proj}_{Y}(\mathrm{Sym}(-))$ (or $\mathbb{P}_{Y}(-):=\mathrm{Projan}_{Y}(\mathrm{Sym}(-))$ for complex analytic case) and composing the natural morphisms $X\cong \mathbb{P}_{X}(\mathcal{L})$ and $\mathbb{P}_{X}(\pi^{*}\pi_{*}\mathcal{L})\cong \mathbb{P}_{Y}(\pi_{*}\mathcal{L})\times_{Y}X\to \mathbb{P}_{Y}(\pi_{*}\mathcal{L})$,
where the operator $\mathbb{P}_{Y}(-)$ on a stack $Y$ can be defined via smooth descent over $Y$. 
A line bundle $\mathcal{L}$ is called {\em very ample} over $Y$ if it is basepoint-free and $\varphi_{\mathcal{L}}$ is a closed immersion.
A line bundle $\mathcal{L}$ on $X$ is called {\em ample} if there exist a smooth covering $\{Y_{i}\to Y\}_{i}$ and positive integers $m_{i}$ such that
 $\mathcal{L}^{\otimes m_{i}}|_{X\times_{Y}Y_{i}}$ is very ample over $Y_{i}$ for each $i$.
 A relative $\R$-line bundle on $X$ over $Y$ is {\em ample} over $Y$ if it is written smooth locally on $Y$ as a positive linear combination of relatively ample line bundles.
 In this section, the morphism $\pi$ is called {\em projective} if there exists an ample $\R$-line bundle $\mathcal{L}$ on $X$ over $Y$.
 Note that any projective morphism is representable by spaces.

For a projective morphism $\pi\colon X\to Y$ of stacks, 
other properties of relative $\R$-line bundles, e.g., {\em nef}, {\em movable}, {\em big}, {\em pseudo-effective} and so on, are defined similarly as ample $\R$-line bundles, which we omit the definitions.
For example, a relative $\R$-line bundle $\mathcal{L}$ on $X$ over $Y$ is nef over $Y$ if for any proper curve $C$ with a morphism $C\to X$ whose image on $Y$ is a point, the degree of the restriction $\mathcal{L}|_{C}$ is non-negative. 
In particular, an $\R$-line bundle $\mathcal{L}$ on $X$ is nef (resp.\ movable, big, pseudo-effective) over $Y$ if and only if there exists a smooth covering $\{Y_{i}\to Y\}_{i}$ from affine schemes or Stein spaces $Y_i$ such that the restriction $\mathcal{L}|_{X\times_{Y}Y_{i}}$ is nef (resp.\ movable, big, pseudo-effective) over $Y_i$ in the usual sense.
\end{defn}

\begin{defn}[Weil divisor]
Let $X$ be a stack which is equidimensional.
An irreducible and reduced Zariski closed substack $Y\subset X$ of codimension one is called a {\em prime divisor} on $X$.
Let $\mathrm{Div}'(X)$ denote the free abelian group generated by prime divisors on $X$.
For a smooth morphism $U\to X$, the pullback of cycles defines a homomorphism $\mathrm{Div}'(X)\to \mathrm{Div}'(U)$.
This correspondance induces a presheaf $U\mapsto \mathrm{Div}'(U)$ on the lisse-\'{e}tale site of $X$.
Its sheafification is denoted by $\mathrm{Div}(-)$.
The global section $\mathrm{Div}(X)=\lim_{U\to X} \mathrm{Div}(U)$ is called {\em the group of (Weil) divisors} on $X$ and its element is called a {\em (Weil) divisor} on $X$.
For $R=\Q$ or $\R$, the notion of $R$-divisors on $X$ can also be defined similarly by using the presheaf $U\mapsto \mathrm{Div}'(U)\otimes_{\Z}R$ instead of $U\mapsto \mathrm{Div}'(U)$.
The group of $R$-divisors on $X$ is denoted by $\mathrm{Div}_{R}(X)$.
\end{defn}

\begin{defn}[Cartier divisor]
Let $X$ be a stack which is reduced and equidimensional.
Let $\mathcal{K}_{X}$ denote the sheaf of rational or meromorphic functions on the lisse-\'{e}tale site of $X$.
A {\em Cartier divisor} $D$ on $X$ means an element of $\mathrm{CDiv}(X):=H^{0}(X,\mathcal{K}_{X}^{*}/\O_{X}^{*})$.
It has a representation 
$$
D=\{(U_{i},f_{i})_{i}\},
$$
where $\{U_{i}\to X\}_{i}$ is a smooth covering from spaces $U_i$ and $f_{i}$ is a non-zero rational or meromorphic function on $U_{i}$ such that $f_{i}/f_{j}$ is an invertible regular function on $U_{i}\times_{X}U_{j}$.
For $R=\Q$ or $\R$, elements of $\mathrm{CDiv_{R}}(X):=H^{0}(X,\mathcal{K}_{X}^{*}/\O_{X}^{*}\otimes_{\Z} R)$ are said to be {\em $R$-Cartier divisors} on $X$.
The {\em cycle map} $\mathrm{CDiv}_{R}(X)\to \mathrm{Div}_{R}(X)$ is defined by
sending $D=\{(U_{i},f_{i})\}$ to the $R$-divisor $D^{W}$ obtained by gluing the principal divisors $\mathrm{div}_{U_{i}}(f_{i})$ on $U_{i}$, which is injective if $X$ is normal.
The short exact sequence
$$
0\to \O^{*}_{X}\otimes_{\Z}R\to \mathcal{K}^{*}_{X}\otimes_{\Z}R \to \mathcal{K}^{*}_{X}/\O^{*}_{X}\otimes_{\Z}R \to 0
$$
induces the homomorphism $\mathrm{CDiv}_{R}(X)\to \mathrm{Pic}_{R}(X)$.
If $R=\Z$, it sends a Cartier divisor $D=\{(U_{i},f_{i})_{i}\}$ to a line bundle $\O_{X}(D)$ defined by gluing $\O_{U_i}f_{i}^{-1}$ on $U_{i}$ naturally.
\end{defn}

\begin{defn}[Ample model]\label{defn--stack-amplemodel}
Let $\pi\colon X\to Y$ be a projective morphism of stacks and assume that $X$ is normal.
Let $D$ be an $\R$-line bundle on $X$ over $Y$.
Let $\varphi\colon X\dasharrow X'$ be a contraction to a normal stack $X'$ over $Y$, that is, there exists a resolution of indeterminacy $X\xleftarrow{p}\widetilde{X}\xrightarrow{q} X'$ of $\varphi$ such that $p$ and $q$ are proper and $q_{*}\O_{\widetilde{X}}=\O_{X'}$.

\smallskip

\noindent
(1)
Suppose that $\varphi$ is birational or bimeromorphic.
Let us take a smooth covering $\{Y_{i}\to Y\}_{i}$ from spaces $Y_{i}$ so that each $D|_{Y_i}$ comes from an $\R$-Cartier divisor on $X\times_{Y}Y_{i}$.
Let $D'_{i}$ denote the pushforward of $D|_{Y_i}$ as $\R$-divisors.
If all the $D'_{i}$ are $\R$-Cartier, then these are glued together to an $\R$-line bundle on $X$ over $Y$, which we denote by $\varphi_{*}D$.
In this case, we say that $\varphi$ has the {\em Cartier pushforward $\varphi_{*}D$ of $D$}.
Note that if $D$ is big (resp.\ pseudo-effective, movable, numerically trivial) over $Y$, then so is the Cartier pushforward $\varphi_{*}D$.

\smallskip

\noindent
(2)
We say that $\varphi$ is an {\em ample model of $D$}
if the following conditions are satisfied:
\begin{itemize}
\item
There exist an ample $\R$-line bundle $D'$ on $X'$ over $Y$ and a resolution of indeterminacy $X\xleftarrow{p} \widetilde{X}\xrightarrow{q} X'$ of $\varphi$ such that 
$$
p^{*}D=q^{*}D'+E
$$
holds as $\R$-line bundles on $\widetilde{X}$ over $Y$ for some effective $\R$-Cartier divisor $E$ on $\widetilde{X}$.
\item
There exists a smooth covering $\{Y_i\to Y\}_{i}$ from affine or Stein spaces $Y_i$ such that each $D|_{Y_{i}}$ are globally $\mathbb{R}$-Cartier.
Moreover, for every $i$ and every member $B$ in the $\R$-linear system $|D|_{Y_i}/Y_{i}|_{\R}$ over $Y_i$, we have $B\ge E|_{Y_i}$.
\end{itemize}
In this case, $X'\times_{Y}Y_{i}$ is an ample model of $D|_{Y_i}$ in the usual sense ({\cite[Definition~3.6.5]{bchm}}).

Note that $\varphi$ is birational or bimeromorphic if and only if $D$ is big over $Y$.
In this case, $D'$ is the Cartier pushforward of $D$, $E$ is exceptional over $X'$ and $\varphi$ is a birational or bimeromorphic contraction
by applying the proof of {\cite[Lemma~3.6.6]{bchm}} smooth locally on $Y$.
Note that for a smooth morphism $U\to Y$,
the base change $X\times_{Y}U\dasharrow X'\times_{Y}U$ is also an ample model of $D|_U$. 
\end{defn}

The following lemmas say that ample models are essentially unique and the existence can be checked smooth locally on $Y$.

\begin{lem} \label{descentdata}
Let $\pi\colon X\to Y$ be a projective morphism of spaces such that $X$ is normal.
Let $D$ be an $\R$-line bundle on $X$.
Then ample models of $D$ are unique up to canonical isomorphisms if exist.
More precisely, for any two ample models $\varphi_{i}\colon X\dasharrow X_{i}$ ($i=1,2$) of $D$, 
 an isomorphism $\sigma_{21}\colon X_{1}\to X_{2}$ over $Y$ can be attached such that $\sigma_{21}\circ \varphi_{1}=\varphi_{2}$ and $\sigma_{21}$ is the identity if $\varphi_{1}=\varphi_{2}$.
Moreover, let $\varphi_{i}\colon X\dasharrow X_{i}$ ($i=1,2,3$) be three ample models of $D$ and $\sigma_{ij}\colon X_{j}\to X_{i}$ be the above isomorphisms.
Then $\sigma_{31}=\sigma_{32}\circ \sigma_{21}$ holds.
\end{lem}

\begin{proof}
The claim follows from the proof of {\cite[Lemma~3.6.6]{bchm}} and a standard argument.
For the convenience of the reader, we give a proof.

Let $\varphi_{i}\colon X\dasharrow X_i$ ($i=1,2$) be two ample models of $D$.
Let $X\xleftarrow{p}\widetilde{X}\xrightarrow{q_i} X_i$ be a common resolution
with $p^{*}D=q_{i}^{*}D'_{i}+E_{i}$, 
where $D'_{i}$ is ample over $Y$ and $E_{i}$ is effective as in the definition of ample models.
Then $E_1=E_2$ holds as in the proof of {\cite[Lemma~3.6.6]{bchm}}.
Let $\overline{X}_{12}$ denote the normalization of the image of $(q_1,q_2)\colon \widetilde{X}\to X_1\times_{Y}X_2$.
Let $q\colon \widetilde{X}\to \overline{X}_{12}$ and $r_{i}\colon \overline{X}_{12}\to X_i$ denote the natural morphisms with $q_{i}=r_{i}\circ q$.
Let $A:=r_{1}^{*}D'_{1}+r_{2}^{*}D'_{2}$.
Then we have
$$
q^{*}A=q_{1}^{*}D'_{1}+q_{2}^{*}D'_{2}=2q_{i}^{*}D'_{i}.
$$
We will show that $r_{1}$ and $r_2$ are isomorphisms.
Assume contrary that $r_1$ is not an isomorphism.
Then there exists an irreducible curve $C$ on $\overline{X}_{12}$ such that $r_{1}(C)$ is a point.
By the definition of $\overline{X}_{12}$, the image $r_{2}(C)$ is a curve.
Hence $r_{1}^{*}D'_{1}\cdot C=0$ and $r_{2}^{*}D'_{2}\cdot C>0$.
In particular, $A\cdot C>0$.
On the other hand, let $\widetilde{C}$ be a curve on $\widetilde{X}$ such that $q(\widetilde{C})=C$.
Then $q^{*}A\cdot C'=2q_{1}^{*}D'_{1}\cdot C'=0$
since $q_{1}$ contracts $C'$ to a point.
Thus we have $A\cdot C=0$, which is a contradiction.
Similarly, $r_{2}$ is also an isomorphism.
We define $\sigma_{21}:=r_{2}\circ r_{1}^{-1}$.
This is independent of the choice of a common resolution $\widetilde{X}$,
and $\sigma_{21}\circ \varphi_{1}=\varphi_{2}$ holds.
By definition, $\sigma_{21}$ is the identity if $\varphi_{1}=\varphi_{2}$.

Let $\varphi_{i}\colon X\dasharrow X_i$ ($i=1,2,3$) be three ample models of $D$
and $X\xleftarrow{p}\widetilde{X}\xrightarrow{q_i} X_i$ be a common resolution
with $p^{*}D=q_{i}^{*}D'_{i}+E_{i}$.
Let $\overline{X}_{123}$ denote the normalization of the image of $\widetilde{X}\to X_{1}\times_{Y}X_{2}\times_{Y}X_{3}$. 
By construction, we have $\overline{X}_{123} \overset{\cong}{\longrightarrow} \overline{X}_{ij}$, where $(i,j)=(1,2), (2,3), (1,3)$. 
Therefore, we have the commutative diagram
$$
\xymatrix@=23pt{
&&X_{1}\ar@/^3pc/[dddrr]^{\sigma_{31}}\ar@/_3pc/[dddll]_{\sigma_{21}}\\
&\overline{X}_{12}\ar[ur]\ar[ddl]&&\overline{X}_{13}\ar[ul]\ar[ddr]\\
&&\overline{X}_{123}\ar[ul]\ar[ur]\ar[d]\\
X_{2}\ar@/_2pc/[rrrr]_{\sigma_{32}}&&\overline{X}_{23}\ar[ll]\ar[rr]&&X_{3}
}
$$
such that all the morphisms are isomorphisms.
The equality $\sigma_{31}=\sigma_{32}\circ \sigma_{21}$ follows by chasing the above commutative diagram. 
\end{proof}

\begin{lem} \label{glueamplemodel}
Let $\pi\colon X\to Y$ be a projective morphism of stacks and assume that $X$ is normal.
Let $D$ be an $\R$-line bundle on $X$ over $Y$.
Then ample models of $D$ are unique up to canonical isomorphisms if exist.
Moreover, if there exists a smooth covering $\{Y_{i}\to Y\}_{i}$ from spaces $Y_i$ such that $X\times_{Y}Y_i$ has an ample model of $D|_{Y_i}$ for each $i$, then $X$ also has an ample model of $D$.
\end{lem}

\begin{proof}
By Lemma~\ref{descentdata},
it suffices to show the existence part of the claim.
Let $\{Y_i\to Y\}_{i}$ be a smooth covering as in the claim
and put $X_i:=X\times_{Y}Y_i$.
We may assume that $D_i:=D|_{Y_i}$ is an $\R$-Cartier divisor on $X_{i}$ by refining the covering $\{Y_i\to Y\}_{i}$ if necessary.
Let $\varphi_{i}\colon X_i \dasharrow X'_{i}$ be an ample model of $D_i$.
Applying Lemma~\ref{descentdata},
the spaces $X'_{i}$ are glued to a stack $X'$ over $Y$.
Indeed, let
$U:=\bigsqcup_{i}Y_{i}$
and $R:=U\times_{Y}U=\bigsqcup_{i,j}Y_{ij}$, where $Y_{ij}:=Y_{i}\times_{Y}Y_{j}$.
Then the disjoint union $X'_{U}:=\bigsqcup_{i}X'_{i}$
is an ample model of $D|_{U}$ 
and 
$$
X'_{R}:=\bigsqcup_{i,j}(X'_{i}\times_{Y_{i}}Y_{ij}),\quad
X''_{R}:=\bigsqcup_{i,j}(X'_{j}\times_{Y_{j}}Y_{ij})
$$
are both ample models of $D|_{R}$.
By Lemma~\ref{descentdata}, there exists a canonical isomorphism $\sigma\colon X'_{R}\cong X''_{R}$
and
the two projections $\mathrm{pr}_{1}\colon X'_{R}\to X'_{U}$ and $\mathrm{pr}_{1}\circ \sigma \colon X'_{R}\to X'_{U}$ form a smooth groupoid and define $X'$ as the quotient stack $[X'_{U}/X'_{R}]$ (cf.\ {\cite[Tag 044O]{Stacks}}).
The contractions $\varphi_{i}\colon X_{i}\dasharrow X'_{i}$ also descend a contraction $\varphi\colon X\dasharrow X'$
by the same argument as above (applying to the graph of $\varphi_{i}$).
Let $D'_{i}$ be the corresponding ample $\R$-line bundle on $X'_{i}$ over $Y_i$.
Since $D'_{i}|_{Y_{ij}}=D'_{j}|_{Y_{ij}}$ as $\R$-line bundles for each $i,j$,
these $D'_{i}$ are glued together to an ample $\R$-line bundle $D'$ on $X'$ over $Y$.
Let $X\xleftarrow{p}\widetilde{X}\xrightarrow{q}X'$ be the resolution of indeterminacy where $\widetilde{X}$ is the normalization of the closure of the graph of $\varphi$.
Let $E:=p^{*}D-q^{*}D'$ as an $\R$-line bundle on $\widetilde{X}$ over $Y$.
Then $E|_{Y_i}$ is represented by an effective $\R$-divisor $E_{i}$ on $\widetilde{X}\times_{Y}Y_{i}$ for each $i$
and the pullbacks of $E_{i}$ and $E_{j}$ to $\widetilde{X}\times_{Y}Y_{ij}$ are equal as $\R$-divisors.
Then $E$ is also represented by an effective $\R$-divisor on $\widetilde{X}$
defined by gluing of $E_{i}$ as $\R$-divisors.
Hence $\varphi\colon X\dasharrow X'$ is an ample model of $D$.
\end{proof}

\begin{defn}[$t$-th output of a $D$-MMP with scaling of $H$] \label{outputMMP}
Let $\pi\colon X\to Y$ be a projective morphism of stacks
and assume that $X$ is normal.
Let $D$ and $H$ be $\R$-line bundles on $X$ over $Y$ such that $H$ is big over $Y$.
Let us assume that the pseudo-effective threshold of $H|_{U}$ with respect to $D|_{U}$, denoted by $\mu$, is constant for any smooth morphism $U\to Y$ from a space $U$ such that $X\times_{Y}U$ is non-empty.
Note that if $X$ is irreducible, then 
this assumption is satisfied 
and $\mu$ equals
the minimal non-negative number such 
that $D+\mu H$ is pseudo-effective over $Y$.
For $t>\mu$, 
a birational or bimeromorphic map $\varphi_{t}\colon X\dasharrow X_{t}$
is called 
a {\em $t$-th output of a $D$-MMP with scaling of $H$ over $Y$}
if it is an ample model of $D+(t-\varepsilon) H$ for sufficiently small $\varepsilon >0$ smooth locally on $Y$, that is, there exist a smooth covering $\{Y_{i}\to Y\}_{i}$ and positive numbers $\{a_{i}\}_{i}$ such that the base change of $\varphi_{t}$ to $Y_{i}$ is an ample model of $D|_{Y_i}+(t-\varepsilon)H|_{Y_i}$ for $0<\varepsilon<a_i$.

We say that the {\em $D$-MMP with scaling of $H$ over $Y$ exists} if there exists a $t$-th output of a $D$-MMP with scaling of $H$ over $Y$ for any $t>\mu$.
\end{defn}

The following statement is a variant of Lemma~\ref{glueamplemodel}.

\begin{lem} \label{glueMMP}
Let $\pi\colon X\to Y$ be a projective morphism of stacks and assume that $X$ is normal.
Let $D$, $H$ and $\mu$ be as in Definition~\ref{outputMMP} and let $t>\mu$.
Then $t$-th outputs of a $D$-MMP with scaling of $H$ over $Y$ are unique up to canonical isomorphisms if exist.
Moreover, if there exists a smooth covering $\{Y_{i}\to Y\}_{i}$ from spaces $Y_i$ such that $X\times_{Y}Y_i$ has 
a $t$-th output of a $D|_{Y_i}$-MMP with scaling of $H|_{Y_i}$ over $Y_i$ for each $i$, then $X$ also has 
a $t$-th output of a $D$-MMP with scaling of $H$ over $Y$.
\end{lem}

\begin{proof}
The uniqueness follows from Lemma~\ref{glueamplemodel}. 
Suppose that there exist a smooth covering $\{Y_{i}\to Y\}_{i}$, positive real numbers $\{a_{i}\}_{i}$, and birational or bimeromorphic maps $\varphi_{t,i}\colon X\dasharrow X_{t,i}$ over $Y_{i}$ such that $\varphi_{t,i}$ is an ample model of $D|_{Y_i}+(t-\varepsilon)H|_{Y_i}$ for any $0<\varepsilon<a_i$.
Then we may apply the proof of  Lemma~\ref{glueamplemodel} to glue these $X_{t,i}$, and we get the $t$-th output $\varphi_{t}\colon X\dasharrow X_{t}$ of a $D$-MMP with scaling of $H$ over $Y$. 
Note that $\varphi_{t}$ is not necessarily an ample model of some $\mathbb{R}$-line bundle because $a_{i}$ depends on $Y_{i}$ and this number can be arbitrary small. 
\end{proof}

Next we consider outputs of $D$-MMP with scaling of $H$ at many times $t$.

\begin{lem}
Let $\pi\colon X\to Y$ be a projective morphism of stacks with $X$ normal.
Let $D$ be a big $\R$-line bundle on $X$ over $Y$ and $H$ be an $\R$-line bundle on $X$ which is nef over $Y$.
Assume that there exist an ample model $\varphi_{1}\colon X\dasharrow X_{1}$ of $D+H$ and an ample model $\varphi_{2}\colon X\dasharrow X_2$ of $D$ with the Cartier pushforward $\varphi_{1*}D$.
Then $\varphi_{2}\circ \varphi_{1}^{-1}\colon X_{1} \dasharrow X_2$ is an ample model of $\varphi_{1*}D$.
\end{lem}

\begin{proof}
By Lemma~\ref{glueamplemodel}, we may assume that $Y$ is a space and both $D$ and $H$ are represented by $\R$-Cartier divisors on $X$.
Let $p\colon \widetilde{X}\to X$ and $p_{i}\colon \widetilde{X}\to X_{i}$ be common resolutions of indeterminacy.
Let us write as $p^{*}(D+H)=p_{1}^{*}(D_1+H_1)+E_1$ and
$p^{*}D=p_{2}^{*}D_{2}+E_2$,
 where $D_i=\varphi_{i*}D$ for $i=1,2$, $H_1=\varphi_{1*}H$ and $E_{i}$ ($i=1,2$) are effective $p_i$-exceptional divisors, respectively.
 Then 
 $$p_{1}^{*}D_1=p_{2}^{*}D_2+(E_2-E_1-F_1)$$
  holds, where $F_1:=p_{1}^{*}H_1-p^{*}H$.
By the negativity lemma, $F_{1}$ and $E_2-E_1-F_1$ are effective.
 Hence $E_2-E_1-F_1$ is $p_2$-exceptional since so is $E_2$ and $E_2-E_1-F_1\le E_2$.
\end{proof}

\begin{cor} \label{MMPfactor}
Let $\pi\colon X\to Y$ be a projective morphism of stacks with $X$ normal.
Let $D$, $H$ and $\mu$ be as in Definition~\ref{outputMMP} and assume that $H$ is nef over $Y$.
Let $t>s>\mu$ and assume that there exist a $t$-th output of a $D$-MMP $\varphi_{t}\colon X\dasharrow X_{t}$ with scaling of $H$ and an $s$-th output of a $D$-MMP $\varphi_{s}\colon X\dasharrow X_s$ with scaling of $H$.
Then $\varphi_{s}\circ \varphi_{t}^{-1}\colon X_{t} \dasharrow X_{s}$ is an $s$-th output of $\varphi_{t*}D$-MMP with scaling of $\varphi_{t*}H$.
\end{cor}

\begin{defn}[Termination of MMP with scaling] 
Let $\pi\colon X\to Y$, $D$, $H$ and $\mu$ be as in Definition~\ref{outputMMP}.
Let us assume that $H$ is nef over $Y$ and
the $D$-MMP with scaling of $H$ over $Y$ exists.
For $t>s>\mu$, the map $\varphi_{s,t}:=\varphi_{s}\circ \varphi_{t}^{-1}\colon X_{t}\dasharrow X_{s}$ is an $s$-th output of a $\varphi_{t*}D$-MMP with scaling of $\varphi_{t*}H$ by Corollary~\ref{MMPfactor}.
Let 
$$
\mathcal{T}(\pi,D,H)\subset \R
$$
denote the subset of real numbers $t$ such that $t>\mu$ and
$\varphi_{t,t+\varepsilon}$ is not an isomorphism for any sufficiently small $\varepsilon>0$.
Note that accumulation points of $\mathcal{T}(\pi,D,H)$ may exist in general.
We say that the $D$-MMP with scaling of $H$ {\em terminates} (resp.\ {\em terminates smooth locally on $Y$}) if $\mathcal{T}(\pi,D,H)$ has no accumulation points in $\mathbb{R}$ (resp.\ there exists a smooth covering $\{Y_{i}\to Y\}_{i}$ such that the $D|_{Y_i}$-MMP with scaling of $H|_{Y_i}$ terminates for each $i$).
If the $D$-MMP with scaling of $H$  terminates, the $t$-th output $X_{t}$ for the minimal number $t$ in $\mathcal{T}(\pi,D,H)$ is called the {\em final output} of the $D$-MMP with scaling of $H$ and denoted by $X_{\mathrm{final}}$.
If the $D$-MMP with scaling of $H$  terminates smooth locally on $Y$, 
then there exist final outputs $(X\times_{Y}Y_{i})_{\mathrm{final}}$
 for some smooth covering $\{Y_{i}\to Y\}_{i}$.
By the same proof of Lemma~\ref{glueamplemodel},
these $(X\times_{Y}Y_{i})_{\mathrm{final}}$ are glued to a stack over $Y$, which is also called the {\em final output} of the $D$-MMP with scaling of $H$ and denoted by $X_{\mathrm{final}}$.
\end{defn}

\begin{rem}
If a $D$-MMP with scaling of $H$ terminates and $D+\lambda H$ is ample over $Y$ for some $\lambda>0$, then $\mathcal{T}(\pi, D, H)$ is a finite set in the open interval $(\mu, \lambda)$.
Let $t_0>t_1>\cdots>t_m$ denote the elements of $\mathcal{T}(\pi, D, H)$.
Then the $D$-MMP with scaling of $H$ can be written as a finite sequence
$$
X=X_{t_0}\dashrightarrow X_{t_1}\dashrightarrow\cdots \dashrightarrow X_{t_m}=X_{\mathrm{final}}.
$$
But the condition that $D+\lambda H$ is ample over $Y$ for some $\lambda>0$ is not satisfied in general even if $H$ is ample over $Y$ without assuming $Y$ is quasi-compact.
In this case, $\mathcal{T}(\pi, D, H)$ has no upper bound.
\end{rem}

\begin{defn}[Singularities of pairs] \label{defn:lcklt}
Let $(X,\Delta)$ be a pair of a stack $X$ and an effective $\R$-divisor $\Delta$ on $X$.
Then the pair is said to be {\em log canonical} (resp.\ {\em Kawamata log terminal}) if 
for any smooth morphism $U\to X$ from an affine or Stein space $U$, the pair $(U, \Delta_{U})$ is log canonical (resp.\ Kawamata log terminal) in the usual sense, where $\Delta_{U}$ is the flat pullback of $\Delta$ as cycles. 
As usual, we call it an {\em lc pair} (resp.~{\em klt pair}) for short. 
In characteristic zero or complex analytic case,
this is equivalent to the condition that the pair $(U,\Delta_{U})$ is log canonical 
(resp.\ Kawamata log terminal) for some smooth covering $U\to X$ from a scheme or complex analytic space $U$.
This follows from the existence of projective log resolutions and the following standard fact: 
For a projective log resolution $\widetilde{X}\to X$ of $(X,\Delta)$ and a smooth morphism $U\to X$, the base change $\widetilde{X}\times _{X}U\to U$ is also a projective log resolution of $(U,\Delta_{U})$.
\end{defn}

\begin{rem}[Log canonical divisors] \label{rem:logcanodiv}
(1) Let $(X,\Delta)$ be an lc pair and $\pi\colon X\to Y$ a projective morphism between stacks.
Then the log canonical divisor $K_{X}+\Delta$ can be regarded as a relative $\R$-line bundle on $X$ over $Y$ as follows:
Let $U\to Y$ be any smooth morphism from a space $U$ such that the base change $(X_{U}, \Delta_{X_U})$ is lc in the usual sense.
In particular, $K_{X_U}+\Delta_{X_U}$ is $\R$-Cartier and hence regarded as an $\R$-line bundle on $X_U$ over $U$.
For any smooth morphism $V\to U$, the difference between $K_{X_V}+\Delta_{X_V}$ and the pullback of $K_{X_U}+\Delta_{X_U}$ on $X_V$ is the dualizing line bundle $K_{X_V/X_U}$, which is the pullback of the dualizing line bundle $K_{V/U}$.
Hence $K_{X_V}+\Delta_{X_V}$ coincides with the pullback of $K_{X_U}+\Delta_{X_U}$ as relative $\R$-line bundles.
Thus the family $\{K_{X_U}+\Delta_{X_U}\}_{U\to X}$ defines a relative $\R$-line bundle on $X$ over $Y$, which is also denoted by $K_X+\Delta$.

\smallskip

\noindent
(2) More generally, the gluing results of the MMP in this section apply to excellent algebraic stacks that admit a dualizing complex.
Suppose that $Y$ is such a stack, equipped with a dualizing complex $\omega_Y^{\bullet}$, meaning that for every smooth morphism $U \to Y$ from a scheme $U$, the pullback $\omega_Y^{\bullet}|_U$ is a dualizing complex in the usual sense.
Let $\pi\colon X\to Y$ be a projective morphism, and let $\Delta$ be an effective $\R$-divisor on $X$.
Although Grothendieck duality for stacks has been developed (cf.\ \cite{neeman}) and canonical divisors on stacks can therefore be defined in the usual way, the singularities of pairs $(X,\Delta)$ and the log canonical divisor $K_X+\Delta$ can be introduced without explicitly using that definition, as follows.
First, observe that the relative dualizing line bundle $\omega_{U/Y}$ for any smooth morphism $U \to Y$ from a scheme $U$ can be defined via smooth descent over $Y$.
Using this, we define the dualizing complex on $U$ by setting 
$$
\omega_{U}^{\bullet}:=\omega_{Y}^{\bullet}|_{U}\otimes \omega_{U/Y}.
$$
Now let $\pi_U \colon X_U \to U$ denote the base change of $\pi$ by $U \to Y$.
The canonical sheaf $\omega_{X_U}$ is then defined as the lowest cohomology sheaf of $\pi_{U}^{!}\omega_{U}^{\bullet}$.
In this way, we can define log canonical (resp.\ Kawamata log terminal) pairs $(X,\Delta)$, exactly as in Definition~\ref{defn:lcklt}.
For a smooth morphism $V\to U$ of schemes over $Y$, we have
$$
\pi_{V}^{!}\omega_{V}^{\bullet}=\pi_{V}^{!}(\omega_{U}^{\bullet}|_{V}\otimes \omega_{V/U})\cong \pi_{V}^{!}(\omega_{U}^{\bullet}|_{V})\otimes \pi_{V}^{*}\omega_{V/U}\cong \varphi^{*}\pi_{U}^{!}\omega_{U}^{\bullet}\otimes \pi_{V}^{*}\omega_{V/U},
$$
where the last isomorphism follows from {\cite[0E9U]{Stacks}}.
Consequently,
$$
K_{X_V}+\Delta_{V}-\varphi^{*}(K_{X_U}+\Delta_U)=K_{X_V}-\varphi^{*}K_{X_U}=\pi_{V}^{*}K_{V/U},
$$
where $\varphi\colon X_V\to X_U$ is the natural projection.
Thus, we may define the log canonical divisor $K_X+\Delta$ as a relative $\mathbb{R}$-line bundle on $X$ over $Y$, in the same manner as in (1).
\end{rem}

\begin{defn}[Minimal model and Mori fiber space]\label{defn--minmodel-stack}
Let $\pi\colon X\to Y$ be a projective morphism of stacks and $(X, \Delta)$ be an lc pair. Then a birational or bimeromorphic map $X\dasharrow X'$ is called a {\em log minimal model} (resp.~a {\em good minimal model}) of $(X, \Delta)$ over $Y$ if it is a log minimal model (resp.~a good minimal model) of $(X, \Delta)$ in the usual sense smooth locally on $Y$, that is, there exists a smooth covering $\{Y_i\to Y\}_{i}$ from affine or Stein spaces $Y_i$ such that the base change $X\times_{Y}Y_i\dasharrow X'\times_{Y}Y_{i}$ is a log minimal model (resp.~a good minimal model) of $(X, \Delta)\times_{Y}Y_i$ over $Y_i$ (around any compact subset of $Y_i$ in the complex analytic case). 

A birational or bimeromorphic map $X\dasharrow X'$ is called a {\em Mori fiber space} of $(X, \Delta)$ over $Y$ if there exists a contraction $X' \to Z$ and a smooth covering $\{Y_i\to Y\}_{i}$ from affine or Stein spaces $Y_i$ such that that the base change $X\times_{Y}Y_i\dasharrow X'\times_{Y}Y_{i}$ satisfies all the conditions of the Mori fiber space except that the relative Picard number is one.
\end{defn}

\begin{rem}\label{rem--MFS-difference}
The above definition of Mori fiber spaces is different from the definition in \cite[Definition 3.1]{eh-mmp} when $Y$ is a Stein space and $\pi\colon X\to Y$ and some $W \subset Y$ satisfies the property (P). 
Even if $Y$ is a point and $(X,\Delta)$ is $\mathbb{Q}$-factorial dlt, the above definition of Mori fiber spaces is different from that in \cite[Definition 2.2]{birkar-flip}. 
However, we adopt this definition for the convenience of the assertion of Theorem~\ref{thm-glue}. 
For the final output of the MMP of Theorem~\ref{thm-glue}, see Remark \ref{rem--output-mmp-stack} below. 
\end{rem}

By applying Lemma~\ref{glueMMP},
the outputs of a $(K_X+\Delta)$-MMP obtained in
Lemma~\ref{lem--mmp-another-termination} (or the corresponding algebraic versions obtained in \cite{hashizumehu})
can be glued together smooth locally on $Y$.
Hence Theorem~\ref{thm-glue} hold:

\begin{proof}[Proof of Theorem~\ref{thm-glue}]
The assertion (1) follows from Lemma~\ref{lem--mmp-another-termination} and Lemma~\ref{glueMMP}.
The assertion (2) follows from Theorem~\ref{thm--strictmmp-intro}, Lemma~\ref{lem--mmp-another-termination} and Lemma~\ref{glueMMP}.
The assertion (3) follows from Theorems~\ref{thm--mmp-birkar-intro}, \ref{thm--mmp-ampleeffective-intro}, Lemma~\ref{lem--mmp-another-termination} and Lemma~\ref{glueMMP}.
\end{proof}

\begin{proof}[Proof of  Corollary \ref{cor--local-criterion-mmp}]
This follows from Lemma \ref{lem--mmp-another-termination} and Lemma \ref{glueMMP}.
\end{proof}

\begin{rem}\label{rem--output-mmp-stack}
As in Definition \ref{defn--minmodel-stack}, the Mori fiber space $X'$ obtained as the final output of the MMP in Theorem~\ref{thm-glue} admits a (not necessarily extremal) log Fano fibration $X'\to Z$ over $Y$.
Indeed, we apply the basepoint-free theorem smooth locally on $Y$ to $K_{X'}+\Delta'+\mu H'$, the Cartier pushforward of $K_{X}+\Delta+\mu H$,
to obtain log Fano fibrations smooth locally on $Y$.
Since the fibration is an ample model of $K_X+\Delta+\mu H$ over $Y$, these fibrations are glued together to the desired fibration by Lemma~\ref{glueamplemodel}.
Similarly, the good minimal model $X'$ obtained as the final output of the MMP in Theorem~\ref{thm-glue} admits an Iitaka fibration $X'\to Z$ over $Y$.
\end{rem}


\end{document}